\documentclass[a4paper,12pt]{article}
\usepackage[cp1251]{inputenc}
\usepackage{amssymb,amsmath,amsthm}
\usepackage[russian]{babel}
\usepackage{indentfirst}
\usepackage{amssymb,amsmath,amsthm}
\usepackage{hyperref,hypbmsec}
\usepackage[pdftex]{graphicx}

\numberwithin{equation}{section}

\newtheorem{Th}{\hskip\parindent Теорема}[section]
\newtheorem{Le}{\hskip\parindent Лемма}[section]

\newtheorem{Zam}{\hskip\parindent Замечание}[section]

\newcommand{\A}{\mathcal{A}}
\newcommand{\T}{\mathcal{T}}
\newcommand{\R}{\mathcal{R}}
\newcommand{\kk}{\emph{k}}

\newcounter{propet}

\renewcommand{\le}{\leqslant}\renewcommand{\ge}{\geqslant}
\renewcommand{\proofname}{Док}

\begin{document}

\author{V.\,A.\,Bykovskii,\quad D.\,A.\,Frolenkov
\footnote{
This work was supported by the Russian Science Foundation under grant [14-11-00335] and performed in the Institute for Applied Mathematics, Far Eastern Branch, Russian Academy of Sciences.}  }
\title{
Asymptotic formulas for the second moments of $L$-series associated to holomorphic cusp forms on the critical line
}
\date{}
\maketitle
\begin{abstract}
New uniform asymptotic formulas are obtained for the second moment of $L$--series of cusp forms of even weight $2k\ge2$  with respect to the congruence subgroup $\Gamma_0(N).$
\\
\noindent

\end{abstract}

\renewcommand{\proofname}{{\bf Доказательство}}

\setcounter{Zam}0
\section{Введение}
Ориентируясь на обозначения и терминологию из монографии \cite{IK}, обозначим через $S_{2k}(N)$ комплексное линейное пространство голоморфных параболических форм целого четного веса $2k\ge2$ относительно конгруэнц-подгруппы $\Gamma_0(N)$. Каждая форма $f$ из $S_{2k}(N)$ раскладывается в ряд Фурье
\begin{equation*}
f(z)=\sum_{n=1}^\infty \rho_f(n)e(nz),
\end{equation*}
где $e(x)=\exp(2\pi ix).$ Ассоциированные с $f$ $L$--ряды
\begin{equation*}
L_f(s)=\sum_{n=1}^\infty \frac{\rho_f(n)}{n^{s+k-1/2}},\qquad
L^{*}_f(s)=\sum_{n=1}^\infty \frac{\overline{\rho_f(n)}}{n^{s+k-1/2}},
\end{equation*}
абсолютно сходящиеся в полосе $\Re s>1,$ определяют голоморфные по $s$ функции на всей плоскости комплексного переменного (см., например, \cite[стр. 378]{IK}).\par
Пространство $S_{2k}(N)$ конечномерно и мы можем выбрать в нем ортонормированный базис $O_{2k}(N)$ относительно скалярного произведения Петерссона
\begin{equation*}
\langle f,g\rangle=\iint\limits_{\Gamma_0(N)\backslash H}f(z)\overline{g(z)}y^{2k-2}dxdy,\qquad
H=\{z=x+iy| x,y\in \mathbb{R}; y>0\}.
\end{equation*}

Изучение свойств $L_f(s)$ на критической прямой $\Re s=1/2$ является одной из важнейших задач теории автоморфных функций (см.\cite{IK}). Следует отметить, что во многих приложениях достаточно оценок в среднем по тем или иным параметрам.
Асимптотические формулы для различных вторых моментов функций $L_f(s)$ необходимы не только для изучения свойств самих L-рядов, но и для решения многих задач аналитической теории чисел. Положим
{\sloppy

}
\begin{equation}\label{W}
W_{2k}(N;t)=\log N+2\gamma-2\log(2\pi)+\frac{\Gamma'(k+it)}{\Gamma(k+it)}+
\frac{\Gamma'(k-it)}{\Gamma(k-it)},
\end{equation}
\begin{equation*}\label{U}
U_{2k}(t)=\frac{\zeta(1+2it)}{(2\pi)^{2it}}\frac{\Gamma(k+it)}{\Gamma(k-it)}+
\frac{\zeta(1-2it)}{(2\pi)^{-2it}}\frac{\Gamma(k-it)}{\Gamma(k+it)},
\end{equation*}
где $\gamma$--постоянная Эйлера, $\Gamma(s)$ и $\zeta(s),$ соответственно, гамма-функция и дзета-функция Римана. Главным результатами статьи являются
\begin{Th}\label{BFtheorem1}
Пусть $k$ и $N$--натуральные числа, $t$--вещественное. Тогда
\begin{gather*}\label{BFTH1}
\frac{\Gamma(2k-1)}{(4\pi)^{2k-1}}\sum_{f\in O_{2k}(N)}\left|L_f(\frac{1}{2}+it)\right|^2=
W_{2k}(N;t)+(-1)^k\delta_{1,N}U_{2k}(t)+V_1(k,N;t)
\end{gather*}
и для любого $\varepsilon>0$
\begin{gather*}\label{BFTH11}
V_1(k,N;t)\ll_{\varepsilon}\frac{1+|t|}{kN}(kN(1+|t|))^{\varepsilon}.
\end{gather*}
\end{Th}
\begin{Th}\label{BFtheorem2}
Пусть $k$ и $N$--натуральные числа, $T$--вещественное число, большее 1. Тогда
\begin{gather*}
\frac{\Gamma(2k-1)}{(4\pi)^{2k-1}}\sum_{f\in O_{2k}(N)}\int_{0}^{T}\left|L_f(\frac{1}{2}+it)\right|^2dt=\\=
\int_{0}^{T}\left(W_{2k}(N;t)+(-1)^k\delta_{1,N}U_{2k}(t)\right)dt+V_2(k,N;T)\label{BFTH2}
\end{gather*}
и для любого $\varepsilon>0$
\begin{gather*}\label{BFTH21}
V_2(k,N;T)\ll_{\varepsilon}\frac{T}{kN}(kNT)^{\varepsilon}.
\end{gather*}
\end{Th}
В формулировках теорем $\delta_{1,N}$--символ Кронекера, равный 1 при $N=1$ и 0 при $N>1.$\par
Эти результаты при $k>1$ вместе с пояснениями по поводу их доказательств опубликованы в~\cite{BF}. Случай $k=1$ там не рассматривался по причине отсутствия в существующих публикациях (известных авторам) формулы свертки для веса 2 (см. \cite{Byk1},
\cite{Byk2}, \cite[теорема 17]{IS}). В настоящей статье этот пробел устраняется и рассматриваются все случаи $k\ge1.$
{\sloppy

}
Различные асимптотические формулы для второго момента рядов $L_f(s)$ или их обобщений при усреднении по базису $O_{2k}(N)$ можно найти в работах \cite{Ak}, \cite{Kam}, \cite{San}. Так А. Санкаранараян при $k>1$ доказал , что
{\sloppy

}
\begin{gather*}\notag
\frac{\Gamma(2k-1)}{(4\pi)^{2k-1}}\sum_{f\in O_{2k}(N)}\int_{0}^{T}\left|L_f(\frac{1}{2}+it)\right|^2dt\ll
NT(\log NT)^3(\log\log T)^2+N^{5+\epsilon}e^{-C\log^2 T}.
\end{gather*}
Таким образом, теорема \ref{BFtheorem2} обобщает и улучшает результат из работы \cite{San}.\par
Из теории Аткина-Ленера следует, что пространство $S_{2k}(N)$ распадается на подпространства новых и старых форм
{\sloppy

}
\begin{equation*}
S_{2k}(N)=S_{2k}^{new}(N)\oplus S_{2k}^{old}(N).
\end{equation*}
Обозначим через $H_{2k}^{*}(N)$ ортогональный базис пространства $S_{2k}^{new}(N).$  Хорошо известно (см.\cite[§ 14.7]{IK}), что при простом $N$ и $k=\{1,2,3,4,5,7\}$ не существует старых форм и, следовательно, $S_{2k}(N)=S_{2k}^{new}(N).$ Соответственно при таких $N$ и $k$ элементы базисов $O_{2k}(N)$  и $H_{2k}^{*}(N)$ отличаются только нормировкой и формулу \eqref{BFTH1} можно рассматривать как частный случай асимптотической формулы для выражения
{\sloppy

}
\begin{gather}\label{second moment primitive}
\frac{\Gamma(2k-1)}{(4\pi)^{2k-1}}\sum_{f\in H^{*}_{2k}(N)}\rho_f(l)\left|L_f(\frac{1}{2}+it)\right|^2,
\end{gather}
где $l$ натуральное. Отметим, что методы доказательства теоремы \ref{BFTH1} переносятся и на случай произвольного $l$. Детальное исследование величин   \eqref{second moment primitive} связано с большим количеством  приложений \cite[§ 26]{IK}. Например, из асимптотических формул для   \eqref{second moment primitive} следуют результаты о пропорции необнуляемости рядов $L_f(1/2).$ Различные результаты об асимптотическом поведении величин \eqref{second moment primitive} или их обобщений можно найти в \cite{Duk}, \cite{KM}, \cite{R2}, \cite{Van}.
В частном случае при $k=1, t=0$ и $N$--простое теорема \ref{BFtheorem1} совпадает с результатом Х. Буя \cite{Bui}. Отметим, что методы данной работы и статьи \cite{Bui} отличаются друг от друга.
{\sloppy

}
\section{Вспомогательные утверждения}
Для натурального $q$ и целого $a$ положим
\begin{gather*}
\delta_q(a)=
\left\{
              \begin{array}{ll}
               1, & \hbox{если $a\equiv0\pmod{q}$} \\
                0, & \hbox{если $a\not\equiv0\pmod{q}$.}
              \end{array}
\right.
\end{gather*}
Для комплексного $s$ и натурального $n$
\begin{equation}\label{tau}
\tau_s(n)=\tau_{-s}(n)=\sum_{n_1n_2=n}\left(\frac{n_1}{n_2}\right)^s=n^s\sigma_{-2s}(n)=n^{-s}\sigma_{2s}(n)
\end{equation}
обобщенная функция числа делителей
\begin{equation*}
\tau(n)=\tau_{0}(n)=\sigma_{0}(n)=\sum_{n_1|n}1.
\end{equation*}
Здесь и в дальнейшем $q,n,n_1$ и $n_2$ натуральные числа.\par
Для целых $m_1,\, m_2$ и $M$
\begin{equation*}
S_q(m_1,m_2;M)=\sum_{a,b=0}^{q-1}\delta_q(ab-M)e\left(\frac{m_1a+m_2b}{q}\right)
\end{equation*}
обобщенная сумма Клостермана. Нам понадобится оценка (см.\cite[следствие 11.12]{IK})
\begin{gather}\notag
\left|S_q(m_1,m_2;\pm1\right|\le \tau(q)(m_1,m_2,q)^{1/2}q^{1/2}\ll_{\varepsilon}
q^{1/2+\varepsilon}(m_1m_2,q)^{1/2}.\label{Kloosterman.bound}
\end{gather}
Для функции Бесселя (см.\cite[глава 7]{BE})
\begin{equation*}
J_s(z)=\sum_{n=0}^{\infty}\frac{(-1)^n}{\Gamma(n+1)\Gamma(n+1+s)}\left(\frac{z}{2}\right)^{s+2n}
\end{equation*}
при $1-2\Re\lambda<\Delta<0$ и положительном вещественном $y$ имеет место интегральное представление (см.\cite[стр. 30, ф. (34)]{BE})
\begin{equation}\label{Barnes}
J_{2\lambda-1}(y)=\frac{1}{4\pi i}\int\limits_{\Re s=\Delta}
\frac{\Gamma(\lambda-1/2+s/2)}{\Gamma(\lambda+1/2-s/2)}\left(\frac{y}{2}\right)^{-s}ds
\end{equation}
и выполняется оценка
\begin{equation*}
J_{2\lambda-1}(y)\ll\min\left\{y^{2\Re\lambda-1},y^{-1/2}\right\},
\end{equation*}
равномерная по $\lambda$ на всех компактых подмножествах в полосе $\Re\lambda>1/2.$  Положим
\begin{equation}\label{K+K-def}
\kk^+(s;z)=\frac{1}{2\sin(\pi s)}\left(J_{-2s}(z)-J_{2s}(z)\right),\qquad
\kk^-(s;z)=\frac{2}{\pi}\cos(\pi s)K_{2s}(z),
\end{equation}
где $K_s(z)$ модифицированная функция Бесселя. При $t>0$ справедливы следующие оценки (см. \cite[ф. (14), (21), (25), (26)]{BookerStrom})
\begin{equation}\label{K.bound>1}
0<K_{it}(tz)\ll e^{-\pi t/2-t\rho(z)}\min\left\{\frac{1}{t^{1/3}},\frac{1}{t^{1/2}(z^2-1)^{1/4}} \right\},
\quad  \hbox{если} \quad z\ge1,
\end{equation}
где $\rho(z)=\sqrt{z^2-1}-\arccos\frac{1}{z}$ и
\begin{equation}\label{K.bound<1}
|K_{it}(tz)|\ll e^{-\pi t/2}\min\left\{\frac{1}{t^{1/3}},\frac{1}{t^{1/2}(1-z^2)^{1/4}} \right\},
\quad  \hbox{если} \quad z<1.
\end{equation}
Также имеет место оценка (см. \cite[ф.(5.16)]{Dunster})
\begin{equation}\label{K+bound}
|\kk^+(it;tz)|\ll\frac{1}{t^{1/2}(1+z^2)^{1/4}}.
\end{equation}
Положим
\begin{equation}\label{phi.gamma.def}
\varphi_s(n)=\sum_{d|n}\frac{\mu(d)}{d^{1+s}},\qquad
\gamma(s;M;q)=\sum_{d|(M,q)}d^{-s}\varphi_s\left(\frac{q}{d}\right).
\end{equation}
Для натурального $M$ и бесконечно дифференцируемой функции $f:(-\infty,\infty)\rightarrow\mathbb{C}$ с компактным носителем, отделенным от нуля, справедлива следующая обобщенная формула суммирования Вороного \cite{BVK}
\begin{gather}\notag
\sum_{n\equiv M(mod\,q)}f(n)\tau_s(n)=
\frac{\zeta(1+2s)}{q}\gamma(2s;M;q)\int_{-\infty}^{\infty}f(x)|x|^sdx+\\\notag+
\frac{\zeta(1-2s)}{q}\gamma(-2s;M;q)\int_{-\infty}^{\infty}f(x)|x|^{-s}dx+\\+
\frac{2\pi}{q^2}\sum_{m,n=1}^{\infty}\left(\frac{m}{n}\right)^{s}\left(
g^{+}\left(16\pi^2\frac{mn}{q^2} \right)S_q(m,n;M)+
g^{-}\left(16\pi^2\frac{mn}{q^2} \right)S_q(m,n;-M)\right),
\label{Voronoi}
\end{gather}
где
\begin{gather*}\label{Voronoitransform+}
g^{\pm}\left(y\right)=\int_{0}^{\infty}\left(
\kk^{\pm}(s;\sqrt{xy})f(x)+\kk^{\mp}(s;\sqrt{xy})f(-x)\right)dx.
\end{gather*}
В дальнейшем нам также понадобится следующее утверждение
\begin{Le}\label{Sum of gamma}
Для $(M,N)=1$ и $\Re s>\max\{1,-\Re v\}$ выполнено следующее равенство
\begin{gather}\label{sumofgamma}
\sum_{q=1}^{\infty}\frac{\gamma(v;M;qN)}{q^{s}}=
\frac{\tau_{(s+v)/2}(M)}{M^{(s+v)/2}}\frac{\zeta(s)}{\zeta(1+s+v)}\frac{\phi_v(N)}{\phi_{s+v}(N)}.
\end{gather}
\begin{proof}
Опираясь на \eqref{phi.gamma.def} и условие $(M,N)=1$, получаем
\begin{gather*}
\sum_{q=1}^{\infty}\frac{\gamma(v;M;qN)}{q^{s}}=
\sum_{q=1}^{\infty}\frac{1}{q^s}\sum_{d|(M,q)}d^{-v}\varphi_v\left(\frac{qN}{d}\right)=
\sum_{d|M}d^{-v-s}\sum_{q=1}^{\infty}\frac{\varphi_v\left(qN\right)}{q^s}=\\=
\sum_{d|M}d^{-v-s}\sum_{q=1}^{\infty}\frac{1}{q^s}\sum_{m|qN}\frac{\mu(m)}{m^{1+v}}=
\sum_{d|M}d^{-v-s}\sum_{m=1}^{\infty}\frac{\mu(m)}{m^{1+v}}\frac{(m,N)^s}{m^s}\zeta(s)=\\=
\frac{\tau_{(s+v)/2}(M)}{M^{(s+v)/2}}\zeta(s)\sum_{m=1}^{\infty}\frac{\mu(m)}{m^{1+v}}\frac{(m,N)^s}{m^s}.
\end{gather*}
Представляя ряд по $m$ в виде эйлеровского произведения, получаем
\begin{gather*}
\sum_{m=1}^{\infty}\frac{\mu(m)}{m^{1+v}}\frac{(m,N)^s}{m^s}=
\prod_{p}\left(1-\frac{(p,N)^s}{p^{1+v+s}}\right)=
\prod_{p\nmid N}\left(1-\frac{1}{p^{1+v+s}}\right)
\prod_{p|N}\left(1-\frac{1}{p^{1+v}}\right)=\\=
\phi_v(N)\prod_{p}\left(1-\frac{1}{p^{1+v+s}}\right)\prod_{p|N}\left(1-\frac{1}{p^{1+v+s}}\right)^{-1}=
\frac{\phi_v(N)}{\zeta(1+s+v)\phi_{s+v}(N)}.
\end{gather*}
Лемма доказана.
\end{proof}
\end{Le}

Дзета-функция Лерха с вещественными параметрами $\alpha$ и $\beta$ в полосе $\Re s>1$ определяется по формуле
{\sloppy

}
\begin{equation*}
\xi(\alpha,\beta;s)=\sum_{n+\alpha>0}\frac{e(n\beta)}{(n+\alpha)^s}
\end{equation*}
с абсолютно сходящимся в полосе $\Re s>1$ рядом из правой части. Она периодична по $\beta$ с периодом 1, продолжается голоморфно на всю плоскость комплексного переменного по $s$ за исключением точки $s=1$ для целого $\beta.$ В последнем случае в этой точке у $\xi(\alpha,\beta;s)$ простой полюс с вычетом 1. Нам понадобятся функции
\begin{equation*}
\xi(\alpha,0;s),\qquad \xi(0,\beta;s),
\end{equation*}
связанные функциональным соотношением
\begin{equation}\label{Lerch.func.equat}
\xi(0,\alpha;1-s)=\frac{\Gamma(s)}{(2\pi)^s}
\left(
e\left(\frac{s}{4}\right)\xi(\alpha,0;s)+e\left(-\frac{s}{4}\right)\xi(-\alpha,0;s)
\right).
\end{equation}
В частном случае $\alpha=0$ мы получаем функциональное уравнение для дзета-функции Римана
\begin{equation}\label{zeta.func.equat}
\zeta(1-s)=\frac{2\Gamma(s)}{(2\pi)^s}
\cos\left(\frac{\pi s}{2}\right)\zeta(s).
\end{equation}
В дальнейшем нам понадобится разложение дзета-функции Римана в ряд Лорана в окрестности точки $s=1$
\begin{equation}\label{zeta.Loran}
\zeta(s)=\frac{1}{s-1}+\gamma+O(|s-1|).
\end{equation}
В частном случае для дзета-функции Лерха при $\Re s>0$ справедливо тождество
\begin{equation*}\label{Lerch.integral.repres.}
\xi\left(0,\frac{a}{q};s\right)=
\sum_{n+a/q>0}^{N-1}\frac{e(na/q)}{(n+a/q)^s}+
s\int_{N}^{\infty}\frac{S(x)}{x^{s+1}}dx,
\end{equation*}
где $q>1,\,(a,q)=1$ и
\begin{equation*}
S(x)=\sum_{n\le x}e\left(\frac{an}{q}\right).
\end{equation*}
Доказательство аналогично \cite[стр.99]{Chudakov} и основывается на применении формулы суммирования Абеля. Выбирая $N=[T],$ получаем  при $0<\sigma<1$ оценку
\begin{equation}\label{Lerch.bound0}
\left|\xi\left(0,\frac{a}{q};\sigma+iT\right)\right|\ll_{\sigma} T^{1-\sigma}.
\end{equation}
Определим величину
\begin{equation*}
\mu(\alpha,\beta;\sigma+it)=\limsup_{t\rightarrow\pm\infty}\frac{\log|\xi(\alpha,\beta;\sigma+it)|}{\log|t|}.
\end{equation*}
В работе \cite{Kat} доказано, что
\begin{equation}\label{Lerchbound}
\mu(\alpha,\beta;\sigma+it)\le\max\left\{\frac{1}{2}-\sigma,\frac{1-\sigma}{2},0\right\}.
\end{equation}

Определим для $\Re s>1+|\Re v|$ абсолютно сходящимся рядом следующую функцию комплексного переменного $s$
\begin{equation*}
G(s,v;q)=\sum_{m,n=1}^{\infty}\frac{S_q(m,n;1)}{(mn)^s}\left(\frac{m}{n}\right)^v.
\end{equation*}
При $q=1$
\begin{equation}\label{G.ifq=1}
G(s,v;1)=\zeta(s-v)\zeta(s+v).
\end{equation}
Из теории дзета-функции Римана следует, что $G(s,v;1)$--голоморфная по $s$ на всей плоскости комплексного переменного функция за исключением точек $s_1=1+v, s_2=1-v$ с полюсами первого порядка при $v\neq0,$ а $\zeta(1+2v)$ и
$\zeta(1-2v)$ соответствующие им вычеты. Функция $G(s,v;1)$ удовлетворяет функциональному уравнению
\begin{equation*}
G(s,v;1)=2\Gamma(1-s+v)\Gamma(1-s-v)(2\pi)^{2s-2}\left(-\cos\pi s+\cos\pi v\right)G(1-s,-v;1).
\end{equation*}
Следующее утверждение является обобщением вышесказанного
\begin{Le}\label{func.equation}
При $q>1$ для любого комплексного $v$ функция $G(s,v;q)$ продолжается по $s$ голоморфно на всю плоскость комплексного переменного и для $\Re s<-|\Re v|$ выполняется равенство
{\sloppy

}

\begin{gather}\notag
G(s,v;q)=2\Gamma(1-s+v)\Gamma(1-s-v)\left(\frac{2\pi}{q}\right)^{2s-2}\\
\times\left(
-\cos\pi s\sum_{m,n=1}^{\infty}\frac{\delta_q(mn-1)}{(mn)^{1-s}}\left(\frac{m}{n}\right)^{-v}
+\cos\pi v\sum_{m,n=1}^{\infty}\frac{\delta_q(mn+1)}{(mn)^{1-s}}\left(\frac{m}{n}\right)^{-v}
\right).\label{G1}
\end{gather}
\begin{proof}
Из определения суммы Клостермана непосредственно следует, что
\begin{equation}\label{Gstart}
G(s,v;q)=\sum_{a,b=0}^{q-1}\delta_q(ab-1)\xi\left(0,\frac{a}{q};s-v\right)\xi\left(0,\frac{b}{q};s+v\right).
\end{equation}
Для $q>1$ выполняется равенство $\delta_q(-1)=0$ и поэтому суммирование по $a$ и $b$ можно вести от 1 до $q-1$. Но в таком случае $\frac{a}{q}$ и $\frac{b}{q}$ не являются целыми числами. Поэтому $G(s,v;q)$ является линейной комбинацией голоморфных функций на всей плоскости комплексного переменного. Так как $\Re s<-|\Re v|$, то используя функциональное уравнение \eqref{Lerch.func.equat} получаем равенство
{\sloppy

}

\begin{gather}\notag
G(s,v;q)=\Gamma(1-s+v)\Gamma(1-s-v)\left(2\pi\right)^{2s-2}
\sum_{a,b=0}^{q-1}\delta_q(ab-1)\\\notag
\times\left(
e\left(\frac{1-s+v}{4}\right)\xi\left(\frac{a}{q},0;1-s+v\right)+
e\left(-\frac{1-s+v}{4}\right)\xi\left(-\frac{a}{q},0;1-s+v\right)\right)\\
\times\left(
e\left(\frac{1-s-v}{4}\right)\xi\left(\frac{b}{q},0;1-s-v\right)+
e\left(-\frac{1-s-v}{4}\right)\xi\left(-\frac{b}{q},0;1-s-v\right)\right).\label{G2}
\end{gather}
Заметим, что
\begin{gather*}
\xi\left(\frac{a}{q},0;s\right)=\sum_{m_1+a/q>0}\frac{1}{(m_1+a/q)^s}=
q^s\sum_{qm_1+a>0}\frac{1}{(qm_1+a)^s}=q^s\sum_{m\equiv a(mod q)}\frac{1}{m^s},
\end{gather*}
\begin{gather*}
\xi\left(\frac{b}{q},0;s\right)=q^s\sum_{n\equiv b(mod q)}\frac{1}{n^s}.
\end{gather*}
Так как для $(m,n)\equiv\pm(a,b)\pmod{q}$
\begin{gather*}
\delta_q(ab-1)=\delta_q((-a)(-b)-1)=\delta_q(mn-1),
\end{gather*}
а для $(m,n)\equiv\pm(-a,b)\pmod{q}$
\begin{gather*}
\delta_q((-a)b-1)=\delta_q(a(-b)-1)=\delta_q(mn+1),
\end{gather*}
то формула \eqref{G2} превращается в \eqref{G1}.
\end{proof}
\end{Le}
Преобразование Меллина от функции $f:[0,\infty)\rightarrow \mathbb{C}$ определяется по формуле
\begin{gather}\label{Mellin.def}
\hat{f}(s)=\int_{0}^{\infty}f(x)x^{s-1}dx.
\end{gather}
Если
\begin{gather*}
f(x)=
\left\{
              \begin{array}{ll}
               O(x^{-a-\varepsilon}), & \hbox{при $x\rightarrow 0+$} \\
               O(x^{-b+\varepsilon}), & \hbox{при $x\rightarrow +\infty$,}
              \end{array}
\right.
\end{gather*}
где $\varepsilon>0$ и $a<b,$ то интеграл в \eqref{Mellin.def} сходится абсолютно и определяет аналитическую функцию в полосе $a<\Re s<b.$ Обратное преобразование Меллина вычисляется по формуле (см., например,\cite{ParisKam})
\begin{gather}\label{Mellin.invers.def}
f(x)=\frac{1}{2\pi i}\int\limits_{\Re s=c}\hat{f}(s)x^{-s}ds,\quad \hbox{где}\quad a<c<b.
\end{gather}
Зафиксируем бесконечно дифференцируемую убывающую функцию $\eta:[0,\infty)\rightarrow[0,1],$ для которой
\begin{gather}\label{eta.def1}
\eta(x)=1,\quad \hbox{если}\quad  0\le x\le \frac{1}{2} \quad\hbox{и}\quad  \eta(x)=0 ,\quad\hbox{если} \quad x\ge 2,
\end{gather}
а также
\begin{gather}\label{eta.def2}
\eta(x)+\eta(1/x)=1.
\end{gather}
\begin{Le}\label{etaMellin}
Преобразование Меллина $\hat{\eta}(s)$ является мероморфной функцией с единственным простым полюсом в точке $s=0$ и для $s\neq0$
{\sloppy

}
\begin{gather}\label{etaMell+etaMell=0}
\hat{\eta}(s)+\hat{\eta}(-s)=0.
\end{gather}
При этом  для любого $\varepsilon>0$ в области
$\min\limits_{j\ge0}dist(s,-j)>\varepsilon$
имеет место оценка
\begin{gather}\label{etaMellbound}
\hat{\eta}(s)=\frac{(-1)^k}{s(s+1)\ldots(s+k-1))}
\int_{1/2}^{2}\eta^{(k)}(x)x^{s+k-1}dx\ll_{\varepsilon}\frac{1}{(1+|s|)^k}.
\end{gather}
\begin{proof}
Заметим, что
\begin{gather*}
\hat{\eta}(s)=\int_{0}^{\infty}\eta(x)x^{s-1}dx=\int_{0}^{1/2}x^{s-1}dx+
\int_{1/2}^{\infty}\eta(x)x^{s-1}dx=
\frac{1}{2^s s}+\int_{1/2}^{\infty}\eta(x)x^{s-1}dx,
\end{gather*}
и последний интеграл задает аналитическую по $s$ функцию. Интегрируя по частям, получим \eqref{etaMellbound}. Чтобы доказать \eqref{etaMell+etaMell=0}, запишем $\hat{\eta}(s)$ и $\hat{\eta}(-s)$ в виде
\begin{gather*}
\hat{\eta}(s)=\frac{-1}{s}\int_{0}^{\infty}\eta'(x)x^{s}dx,\\
\hat{\eta}(-s)=\frac{1}{s}\int_{0}^{\infty}\eta'(x)x^{-s}dx=
\frac{1}{s}\int_{0}^{\infty}\eta'\left(\frac{1}{x}\right)x^{s-2}dx.
\end{gather*}
Отсюда находим, что
\begin{gather*}
\hat{\eta}(s)+\hat{\eta}(-s)=\frac{-1}{s}\int_{0}^{\infty}
\left(\eta'(x)-\frac{1}{x^2}\eta'\left(\frac{1}{x}\right)\right)x^{s}dx.
\end{gather*}
Дифференцируя \eqref{eta.def2} мы получаем, что $\eta'(x)-\frac{1}{x^2}\eta'\left(\frac{1}{x}\right)=0$ и, следовательно, соотношение \eqref{etaMell+etaMell=0} доказано.
{\sloppy

}
\end{proof}
\end{Le}

\section{Формула свертки}
Пусть $u$ и $v$ комплексные числа с $|\Re u|<k-1$ и $\Re v=0.$ Формула свертки (см.\cite{Byk1} или \cite{Byk2}) имеет следующий вид
\begin{gather}\label{convolution k>2}
\frac{\Gamma(2k-1)}{(4\pi)^{2k-1}}\sum_{f\in O_{2k}(N)}L_f(\frac{1}{2}+u+v)L^{*}_f(\frac{1}{2}+u-v)=\\\notag=
\zeta(1+2u)+
\left(\frac{2\pi}{\sqrt{N}}\right)^{4u}\frac{\Gamma(k-u+v)\Gamma(k-u-v)}{\Gamma(k+u+v)\Gamma(k+u-v)}\zeta(1-2u)+\\\notag+
(-1)^k\delta_{1,N}\left(
\left(\frac{2\pi}{\sqrt{N}}\right)^{2u-2v}\frac{\Gamma(k-u+v)}{\Gamma(k+u-v)}\zeta(1+2v)+
\left(\frac{2\pi}{\sqrt{N}}\right)^{2u+2v}\frac{\Gamma(k-u-v)}{\Gamma(k+u+v)}\zeta(1-2v)
\right)+\\\notag+
2(2\pi)^{2u}\cos\pi u\sum\limits_{m_1m_2-n_1n_2=1\atop n_1\equiv0(mod N)}\left(\frac{m_1m_2}{n_1n_2}\right)^{u}
\left(\frac{n_2}{n_1}\right)^{u}\left(\frac{m_2}{m_1}\right)^{v}
\frac{1}{(m_1m_2)^k}H_k\left(u,v;\frac{1}{m_1m_2}\right)+\\\notag+
2(2\pi)^{2u}(-1)^k\cos\pi v\sum\limits_{m_1m_2-n_1n_2=-1\atop n_1\equiv0(mod N)}\left(\frac{m_1m_2}{n_1n_2}\right)^{u}
\left(\frac{n_2}{n_1}\right)^{u}\left(\frac{m_2}{m_1}\right)^{v}
\frac{1}{(m_1m_2)^k}H_k\left(u,v;\frac{-1}{m_1m_2}\right),
\end{gather}
где
\begin{gather}\label{H.func}
H_k\left(u,v;y\right)=\frac{\Gamma(k-u+v)\Gamma(k-u-v)}{\Gamma(2k)}{}_2F_{1}(k-u+v,k-u-v,2k;y)
\end{gather}
с гипергеометрической функцией Гаусса $_2F_{1}(a,b,c;y),$ а участвующие в суммированиях переменные $m_1,m_2,n_1,n_2$--натуральные числа.\par
Ограничение $|\Re u|<k-1$ предполагает, что $k\ge2,$ и именно для таких $k$ доказана формула свертки в  \cite{Byk2} и \cite{Byk1}.
В случае $k=1$ (вес 2) ситуация более сложная. Это связано с тем, что при $y\rightarrow0$
\begin{gather}\label{H1.asym.for}
\frac{\Gamma(2)}{\Gamma(1-u+v)\Gamma(1-u-v)}H_1\left(u,v;y\right)=1+O(y),
\end{gather}
и, следовательно, каждое из двух последних слагаемых расходящиеся ряды.  С данной проблемой мы справимся с помощью замены натурального $k$ комплексной переменной $\lambda$ и дальнейшим аналитическим продолжением в точку $\lambda=k=1.$ \par
Доказательство формулы свертки базируется на классической формуле следа Петерссона (см.\cite[утв. 14.5]{IK})
{\sloppy

}
\begin{equation*}\label{Petersson}
\frac{\Gamma(2k-1)}{(4\pi)^{2k-1}}\sum_{f\in O_{2k}(N)}\frac{\overline{\rho_f(m)}\rho_f(n)}{(mn)^{k-1/2}}=
\delta_{m,n}+2\pi(-1)^k\sum_{q\equiv0(mod N)}\frac{1}{q}S_q(m,n;1)
J_{2k-1}\left(4\pi\frac{\sqrt{mn}}{q}\right),
\end{equation*}
где $m$ и $n$ любые натуральные числа.\par
Пусть $u$ и $v$ такие, что $\Re u>3/4,\,\Re v=0.$ Действуя формально, умножим обе части формулы следа Петерссона на
\begin{gather*}
\frac{1}{(mn)^{1/2+u}}\left(\frac{m}{n}\right)^v,
\end{gather*}
а затем просуммируем по всем натуральным $m$ и $n.$ В результате получим равенство
\begin{gather}\label{basicformula1}
\frac{\Gamma(2k-1)}{(4\pi)^{2k-1}}\sum_{f\in O_{2k}(N)}
L_f(\frac{1}{2}+u+v)L^{*}_f(\frac{1}{2}+u-v)=
\zeta(1+2u)+2\pi(-1)^kD_N(u,v;k),
\end{gather}
где
\begin{gather}\label{D0function}
D_N(u,v;\lambda)=
\sum_{m,n=1}^{\infty}\frac{1}{(mn)^{1/2+u}}\left(\frac{m}{n}\right)^v
\sum_{q\equiv0(mod N)}\frac{S_q(m,n;1)}{q}J_{2\lambda-1}\left(4\pi\frac{\sqrt{mn}}{q}\right).
\end{gather}

\begin{Zam}\label{Zam.lambda=k}
В определении $D_N(u,v;\lambda)$ мы вместо натурального $k$ рассматриваем комплексный параметр $\lambda.$ Это позволит нам получить формулу свертки для веса 2 путем аналитического продолжения по $\lambda$ в точку $\lambda=k=1.$
{\sloppy

}
\end{Zam}
\begin{Le}\label{holom.continuation}
Ряд, определяющий в \eqref{D0function} функцию $D_N(u,v;\lambda)$, абсолютно сходится для
\begin{gather}\label{oblasti1}
\Re \lambda>\frac{3}{4},\qquad \Re u> \frac{3}{4}
\end{gather}
и определяет в этой области голоморфную по $u,v$ и $\lambda$ функцию $D_N(u,v;\lambda).$
\begin{proof}
Пусть $\delta$ любое фиксированное положительное число из интервала $(0,1/4)$ и
\begin{gather}\label{oblasti2}
\delta^{-1}>\Re \lambda>\frac{3}{4}+\delta,\qquad \Re u> \frac{3}{4}+\delta.
\end{gather}
Из оценок для сумм Клостермана и функции Бесселя следует, что интересующий нас ряд мажорируется следующим рядом с вещественными положительными слагаемыми (и некоторой константой $c(\delta)>0$)
\begin{gather*}
c(\delta)\sum_{m,n=1}^{\infty}\frac{1}{(mn)^{1/2+\Re u}}
\sum_{q=1}^{\infty}\frac{\tau(q)(mn,q)^{1/2}}{q^{1/2}}
\min\left\{\left(\frac{\sqrt{mn}}{q}\right)^{2\Re \lambda-1},
\left(\frac{q}{\sqrt{mn}}\right)^{1/2}\right\}.
\end{gather*}
Объединим слагаемые $m$ и $n$ с одинаковым произведением $mn=d.$ По этой причине последний ряд оценивается величиной
\begin{gather*}
c(\delta)\sum_{d=1}^{\infty}\frac{\tau(d)}{d^{3/2+\delta}}
\sum_{q<\sqrt{d}}\tau(q)(d,q)^{1/2}+
c(\delta)\sum_{d=1}^{\infty}\frac{\tau(d)}{d^{7/4+\delta-\Re \lambda}}
\sum_{q\ge\sqrt{d}}\frac{\tau(q)(d,q)^{1/2}}{q^{2\Re \lambda-1/2}}\ll_{\delta}\\\ll_{\delta}
\sum_{d=1}^{\infty}\frac{\tau(d)}{d^{1+\delta/2}}\ll_{\delta}1.
\end{gather*}
Из полученной оценки следует абсолютная сходимость ряда, определяющего функцию $D_N(u,v;\lambda)$ и ее голоморфность по всем трем комплексным переменным $u,v$ и $\lambda$ в области \eqref{oblasti2}. Так как $\delta$ любое положительное число из интервала $(0,1/4),$ то утверждение леммы \ref{holom.continuation} выполняется в полном объеме.
\end{proof}
\end{Le}
Ввиду абсолютной сходимости  мы можем  в области \eqref{oblasti1} переставить порядки суммирования и записать \eqref{D0function} в виде
{\sloppy

}

\begin{gather}\label{D0func.transform}
D_N(u,v;\lambda)=\sum_{q\equiv0(mod N)}\frac{1}{q}D_N(u,v,\lambda;q),
\end{gather}
где
\begin{gather*}\label{D1function1}
D_N(u,v,\lambda;q)=
\sum_{m,n=1}^{\infty}\frac{1}{(mn)^{1/2+u}}\left(\frac{m}{n}\right)^v
S_q(m,n;1)J_{2\lambda-1}\left(4\pi\frac{\sqrt{mn}}{q}\right).
\end{gather*}
\begin{Le}\label{holom.continuation}
Для $\Re \lambda>3/4, \Re u>3/4$ и $\Re v=0$
\begin{gather}\label{D1func.int.repr.1}
D_N(u,v,\lambda;q)=
\frac{1}{4\pi i}\int\limits_{\Re s=\Delta}
\frac{\Gamma(\lambda-1/2+s/2)}{\Gamma(\lambda+1/2-s/2)}
G\left(\frac{1}{2}+u+\frac{s}{2},v;q\right)\left(\frac{q}{2\pi}\right)^{s}ds,
\end{gather}
при условии $\max\{1-2\Re\lambda,1-2\Re u\}<\Delta<0.$
\begin{proof}
Формула \eqref{D1func.int.repr.1} непосредственно следует из интегрального представления \eqref{Barnes} для функции Бесселя и стандартных оценок (далее $t$ - вещественное)
{\sloppy

}
\begin{gather}\label{Gammabound}
\left|\frac{\Gamma(\lambda-1/2+\Delta/2+it/2)}{\Gamma(\lambda+1/2-\Delta/2-it/2)}\right|\ll_{\lambda,\Delta}
\frac{1}{(1+|t|)^{1-\Delta}},
\end{gather}
\begin{gather}\label{Gbound}
G\left(\frac{1}{2}+u+\frac{\Delta}{2}+\frac{it}{2},v;q\right)\ll_{q,u,v}
\sum_{m,n=1}^{\infty}\frac{1}{(mn)^{1/2+\Delta/2+\Re u}},
\end{gather}
равномерных по $\lambda$ на компактных подмножествах. Условие $1-2\Re u<\Delta$ необходимо для абсолютной сходимости ряда в правой части \eqref{Gbound}.
\end{proof}
\end{Le}
\begin{Le}\label{holom.continuation2}
Для $\Re \lambda-1>\Re u>3/4$ и $\Re v=0$
\begin{gather*}\notag
D_N(u,v,\lambda;q)=
\frac{1}{4\pi i}\int\limits_{\Re s=\Delta}
\frac{\Gamma(\lambda-1/2+s/2)}{\Gamma(\lambda+1/2-s/2)}
G\left(\frac{1}{2}+u+\frac{s}{2},v;q\right)\left(\frac{q}{2\pi}\right)^{s}ds+\\+
\delta_{1,q}\left(
\frac{\zeta(1+2v)}{\left(2\pi\right)^{1-2u+2v}}\frac{\Gamma(\lambda-u+v)}{\Gamma(\lambda+u-v)}+
\frac{\zeta(1-2v)}{\left(2\pi\right)^{1-2u-2v}}\frac{\Gamma(\lambda-u-v)}{\Gamma(\lambda+u+v)},
\right)\label{D1func.int.repr.2}
\end{gather*}
при условии $1-2\Re\lambda<\Delta<-1-2\Re u.$
\begin{proof}
В интеграле из \eqref{D1func.int.repr.1} сдвинем прямую интегрирования в область $$1-2\Re\lambda<\Delta<-1-2\Re u.$$ Из \eqref{G.ifq=1} и леммы \ref{func.equation} следует, что при этом мы пройдем два полюса в точках $s_{\pm}=1-2u\pm2v,$  возникающих только при $q=1.$ Обоснуем возможность сдвига прямой интегрирования. Покажем, что можно сдвинуть контур интегрирования с прямой $\Re s=\Delta_1=1-2\Re u+\varepsilon$ на прямую $\Re s=\Delta_2=-1-2\Re u-\varepsilon,$ где $\varepsilon>0.$ Выберем большой параметр $T>0.$
Из представления \eqref{Gstart} и оценок \eqref{Lerchbound}, \eqref{Gammabound} следует, что интегралы по вертикальным лучам $\Re s=\Delta_{1,2}, |\Im s|>T$ сходятся абсолютно и оцениваются величиной
$O\left(T^{1-2\Re u+\epsilon}\right).$ Таким образом нам осталось показать, что интегралы по горизонтальным отрезкам
{\sloppy

}
\begin{gather*}\label{horizontal integrals}
\int_{\Delta_2\pm iT}^{\Delta_1\pm iT}
\frac{\Gamma(\lambda-1/2+s/2)}{\Gamma(\lambda+1/2-s/2)}
G\left(\frac{1}{2}+u+\frac{s}{2},v;q\right)\left(\frac{q}{2\pi}\right)^{s}ds
\end{gather*}
стремятся к нулю при $T\rightarrow\infty.$ После замены переменных получим интегралы
\begin{gather*}
\int_{-\varepsilon\pm iT_1}^{1+\varepsilon\pm iT_1}
\frac{\Gamma(\lambda-u-1+z)}{\Gamma(\lambda+u+1-z)}
G\left(z,v;q\right)\left(\frac{q}{2\pi}\right)^{2z}dz
\end{gather*}
с $T_1=T/2+\Im u,$ которые с помощью формулы Стирлинга оцениваются величиной
\begin{gather}\label{horizontal integrals2}
T_1^{-2-2\Re u}\int_{-\varepsilon}^{1+\varepsilon}
T_1^{2x}\left|G\left(x+iT_1,v;q\right)\right|dx
\end{gather}
Разобьем интеграл \eqref{horizontal integrals2} на два $I_1$ и $I_2$ в соответствии с условиями $x\le\varepsilon$ и $x\ge\varepsilon.$
Используя \eqref{Lerch.bound0}, получаем оценку
$I_2\ll T_1^{-2\Re u+2\varepsilon}.$
Чтобы оценить интеграл $I_1$ мы воспользуемся формулой \eqref{G2}. Возникшие дзета-функции Лерха запишем в виде \cite[лемма 3, стр.24]{VoroninKar}
\begin{gather*}
\xi\left(\alpha,0;s\right)=\sum_{n=0}^{N}\frac{1}{(n+\alpha)^s}+\frac{1}{s-1}
\left(N+\frac{1}{2}+\alpha\right)^{1-s}+s\int_{N+1/2}^{\infty}\frac{1/2-\{u\}}{(u+\alpha)^{s+1}}du,
\end{gather*}
где $\{u\}$ -- дробная часть числа $u.$ Выбирая $N=T_1$ и оценивая получившееся выражение по модулю, получаем оценку
$I_1\ll T_1^{-1-2\Re u+\varepsilon}.$
Таким образом, если $\Re u>3/4,$ то оба интеграла $I_1$ и $I_2$ стремятся к нулю при $T\rightarrow\infty.$ Следовательно, проведенный сдвиг прямой интегрирования законен.
\end{proof}
\end{Le}
Положим
\begin{gather}\label{d.u.v.lam}
d(u,v,\lambda)=\frac{\zeta(1+2v)}{\left(2\pi\right)^{1-2u+2v}}\frac{\Gamma(\lambda-u+v)}{\Gamma(\lambda+u-v)}+
\frac{\zeta(1-2v)}{\left(2\pi\right)^{1-2u-2v}}\frac{\Gamma(\lambda-u-v)}{\Gamma(\lambda+u+v)},
\end{gather}
\begin{gather*}\label{Gamma.u.v.lam}
\Gamma(u,v,\lambda;s)=\frac{\Gamma(\lambda-1/2+s/2)}{\Gamma(\lambda+1/2-s/2)}\Gamma(1/2-u+v-s/2)\Gamma(1/2-u-v-s/2).
\end{gather*}
Используя формулу Стирлинга, при $\Re v=0$ получаем оценку
\begin{gather}\label{Gamma.u.v.lam.bound}
|\Gamma(u,v,\lambda;\sigma+it)|\ll\frac{e^{-\pi|t|/2}}{|t|^{1+2\Re u}}.
\end{gather}
\begin{Le}\label{holom.continuation3}
Для $\Re \lambda-1>\Re u>3/4$ и $\Re v=0$
\begin{gather}\notag
D_N(u,v;\lambda)=
\frac{(2\pi)^{2u-1}}{2\pi i}\int\limits_{\Re s=\Delta}\Gamma(u,v,\lambda;s)\\\notag
\times\Biggl(
\sin\pi\left(u+\frac{s}{2}\right)\sum_{q\equiv0(modN)}\sum_{m,n=1}^{\infty}
\frac{\delta_q(mn-1)}{q^{2u}(mn)^{1/2-u-s/2}}\left(\frac{m}{n}\right)^{-v}+\\
+\cos\pi v\sum_{q\equiv0(mod N)}\sum_{m,n=1}^{\infty}
\frac{\delta_q(mn+1)}{q^{2u}(mn)^{1/2-u-s/2}}\left(\frac{m}{n}\right)^{-v}
\Biggr)ds+
\delta_{1,N}d(u,v,\lambda),\label{D0func.int.repr.1}
\end{gather}
при условии $1-2\Re\lambda<\Delta<-1-2\Re u.$
\begin{proof}
Подставим в формулу \eqref{D0func.transform} результат  леммы \ref{holom.continuation2} и леммы \ref{func.equation}.  Ввиду условий леммы и оценки \eqref{Gamma.u.v.lam.bound} получившиеся суммы и интеграл абсолютно сходятся. Это обеспечивает возможность перестановки суммирования по $q$ и интегрирования по $s.$
\end{proof}
\end{Le}
\begin{Le}\label{holom.continuation4}
Для $\Re \lambda-1>\Re u>3/4$ и $\Re v=0$
\begin{gather}\notag
D_N(u,v;\lambda)=\delta_{1,N}d(u,v,\lambda)+
\frac{(2\pi)^{2u-1}}{2\pi iN^{2u}}\int\limits_{\Re s=\Delta}\Gamma(u,v,\lambda;s)\Biggl(
\zeta(2u)\sin\pi\left(u+\frac{s}{2}\right)+\\
+\sin\pi\left(u+\frac{s}{2}\right)
\sum_{d=1}^{\infty}\frac{\tau_u(d)\tau_v(dN+1)}{d^{u}(dN+1)^{1/2-u-s/2}}
+\cos\pi v\sum_{dN>1}^{\infty}\frac{\tau_u(d)\tau_v(dN-1)}{d^{u}(dN-1)^{1/2-u-s/2}}
\Biggr)ds,\label{D0func.int.repr.2}
\end{gather}
при условии $1-2\Re\lambda<\Delta<-1-2\Re u.$
\begin{proof}
В первой тройной сумме из  \eqref{D0func.int.repr.1} выделяем слагаемое с $m=n=1.$ В оставшейся сумме делаем замену $mn=1+dN$ и, используя \eqref{tau}, получаем
\begin{gather*}
\sum_{q\equiv0(modN)}\sum_{m,n=1}^{\infty}
\frac{\delta_q(mn-1)}{q^{2u}(mn)^{1/2-u-s/2}}\left(\frac{m}{n}\right)^{-v}=
\frac{\zeta(2u)}{N^{2u}}+\\+
\sum_{d=1}^{\infty}\frac{1}{(dN+1)^{1/2-u-s/2}}
\sum_{q\equiv0(mod N)}\frac{\delta_q(dN)}{q^{2u}}\sum_{mn=1+dN}\left(\frac{m}{n}\right)^{-v}=\\=
\frac{\zeta(2u)}{N^{2u}}+
\frac{1}{N^{2u}}
\sum_{d=1}^{\infty}\frac{\tau_u(d)\tau_v(dN+1)}{d^{u}(dN+1)^{1/2-u-s/2}}
\end{gather*}
Во второй тройной сумме из  \eqref{D0func.int.repr.2} делаем замену $mn=-1+dN.$ Заметим, что так как $mn\ge1,$ то при $N=1$ выполнено неравенство $d\ge2.$ Далее действуем аналогично.
\end{proof}
\end{Le}

\begin{Le}\label{holom.continuation5}
Для $\Re \lambda>1,\, \Re u=0,\,\Re v=0$ и $u\neq0$
\begin{gather}\notag
D_N(u,v;\lambda)=
\frac{2(2\pi)^{2u-1}}{N^{2u}}\zeta(2u)\Gamma(2u)\cos\pi(\lambda-u)
\frac{\Gamma(\lambda-u+v)\Gamma(\lambda-u-v)}{\Gamma(\lambda+u+v)\Gamma(\lambda+u-v)}
+\\\notag
+\frac{2(2\pi)^{2u-1}}{N^{2u}}
\Biggl(
\cos\pi(\lambda-u)\sum_{d=1}^{\infty}
\frac{\tau_u(d)\tau_v(dN+1)}{d^{u}(dN+1)^{\lambda-u}}H_{\lambda}\left(u,v;\frac{1}{dN+1}\right)+\\
+\cos\pi v\sum_{d=1+\delta_{1,N}}^{\infty}
\frac{\tau_u(d)\tau_v(dN-1)}{d^{u}(dN-1)^{\lambda-u}}H_{\lambda}\left(u,v;\frac{-1}{dN-1}\right)
\Biggr)+
\delta_{1,N}d(u,v,\lambda),\label{D0func.int.repr.3}
\end{gather}
где функция $H_{\lambda}(u,v;y)$ определена в \eqref{H.func}.
\begin{proof}
Сначала докажем формулу \eqref{D0func.int.repr.3} для $\Re \lambda-1>\Re u>3/4$ и $\Re v=0.$ Для этого достаточно  вычислить следующие интегралы из \eqref{D0func.int.repr.2}
\begin{gather*}
I_1=\frac{1}{2\pi i}\int\limits_{\Re s=\Delta}\Gamma(u,v,\lambda;s)\sin\pi\left(u+\frac{s}{2}\right)ds,\\
I_2=\frac{1}{2\pi i}\int\limits_{\Re s=\Delta}\Gamma(u,v,\lambda;s)\sin\pi\left(u+\frac{s}{2}\right)z^{s/2}_+ds,\\
I_3=\frac{1}{2\pi i}\int\limits_{\Re s=\Delta}\Gamma(u,v,\lambda;s)z^{s/2}_-ds,
\end{gather*}
где $1-2\Re\lambda<\Delta<-1-2\Re u$ и $z_{\pm}=dN\pm1.$ Так как $z_{\pm}\ge1,$ то вычисляя вычеты подынтегральных функций в точках $s_j=1-2\lambda-2j,\, j=0,1,2\ldots,$ получаем
\begin{gather*}
I_2=\frac{2\cos\pi(\lambda-u)}{z_{+}^{\lambda-1/2}}
\frac{\Gamma(\lambda-u+v)\Gamma(\lambda-u-v)}{\Gamma(2\lambda)}{}_2F_{1}\left(\lambda-u+v,\lambda-u-v,2\lambda;\frac{1}{z_{+}}\right),\\
I_3=\frac{2}{z_{-}^{\lambda-1/2}}
\frac{\Gamma(\lambda-u+v)\Gamma(\lambda-u-v)}{\Gamma(2\lambda)}{}_2F_{1}\left(\lambda-u+v,\lambda-u-v,2\lambda;\frac{-1}{z_{-}}\right).
\end{gather*}
Используя равенство
\begin{gather*}
_2F_{1}\left(a,b,c;1\right)=\frac{\Gamma(c)\Gamma(c-a-b)}{\Gamma(c-a)\Gamma(c-b)},
\end{gather*}
получаем, что
\begin{gather*}
I_1=2\cos\pi(\lambda-u)\Gamma(2u)
\frac{\Gamma(\lambda-u+v)\Gamma(\lambda-u-v)}{\Gamma(\lambda+u+v)\Gamma(\lambda+u-v)}.
\end{gather*}
Тем самым формула \eqref{D0func.int.repr.3} для $\Re \lambda-1>\Re u>3/4$ и $\Re v=0$ доказана. Далее по принципу аналитического продолжения получаем утверждение леммы.
\end{proof}
\end{Le}
\begin{Zam}
Применяя в формуле \eqref{D0func.int.repr.3} функциональное соотношение для дзета-функции Римана и полагая $\lambda=k>1,$ мы, используя равенство \eqref{basicformula1}, получим \eqref{convolution k>2}. Однако в случае $k=1$ мы не можем непосредственно положить $\lambda=1$ в формуле \eqref{D0func.int.repr.3}. Получению аналитического продолжения в точку $\lambda=1$ будет посвящен следующий параграф.
{\sloppy

}
\end{Zam}
\section{Аналитическое продолжение в точку $\lambda=1$}
Пусть $P=(1+|t|)^aN^b$ с некоторыми положительными $a$ и $b,$ значения которых будут определены в дальнейшем. Зафиксируем две бесконечно дифференцируемые функции $\alpha,\beta:[0,\infty)\rightarrow[0,1],$ для которых
{\sloppy

}
\begin{gather*}
\beta(x)=1,\quad \hbox{если}\quad  0\le x\le P \quad\hbox{и} \quad \beta(x)=0 ,\quad\hbox{если} \quad x\ge 2P,\\
\alpha(x)=0,\quad \hbox{если}\quad  0\le x\le P \quad\hbox{и}\quad  \alpha(x)=1,\quad \hbox{если} \quad x\ge 2P,\\
\beta(x)+\alpha(x)=1.
\end{gather*}
Подставив данное разбиение единицы в оба ряда из \eqref{D0func.int.repr.3}, получим
\begin{gather}\notag
D_N(u,v;\lambda)=
\frac{2(2\pi)^{2u-1}}{N^{2u}}\zeta(2u)\Gamma(2u)\cos\pi(\lambda-u)
\frac{\Gamma(\lambda-u+v)\Gamma(\lambda-u-v)}{\Gamma(\lambda+u+v)\Gamma(\lambda+u-v)}
+\\+
S_{\beta}(u,v;\lambda)+S_{\alpha}(u,v;\lambda)+
\delta_{1,N}d(u,v,\lambda),\label{D0func.int.repr.4}
\end{gather}
где для $f(x)=\alpha(x)$ или $f(x)=\beta(x)$
\begin{gather}\notag
S_{f}(u,v;\lambda)=
\frac{2(2\pi)^{2u-1}}{N^{2u}}
\Biggl(
\cos\pi(\lambda-u)\sum_{d=1}^{\infty}
\frac{\tau_u(d)\tau_v(dN+1)}{d^{u}(dN+1)^{\lambda-u}}H_{\lambda}\left(u,v;\frac{1}{dN+1}\right)f(dN+1)+\\
+\cos\pi v\sum_{d=1+\delta_{1,N}}^{\infty}
\frac{\tau_u(d)\tau_v(dN-1)}{d^{u}(dN-1)^{\lambda-u}}H_{\lambda}\left(u,v;\frac{-1}{dN-1}\right)f(dN-1)
\Biggr).\label{S.alpha.beta.def}
\end{gather}
Основная и единственная трудность, это аналитическое продолжение функции $S_{\alpha}(u,v;\lambda)$ в точку $\lambda=1.$ На первом шаге мы выделяем "главные части"\, в  $S_{\alpha}(u,v;\lambda),$ а именно
\begin{gather}\notag
S_{\alpha}(u,v;\lambda)=
\frac{2(2\pi)^{2u-1}}{N^{2u}}
\left(
\cos\pi(\lambda-u)C^+\left(u,v,\lambda\right)+
\cos\pi v C^-\left(u,v,\lambda;\right)\right)+\\+
\frac{2(2\pi)^{2u-1}}{N^{2u}}
\frac{\Gamma(\lambda-u+v)\Gamma(\lambda-u-v)}{\Gamma(2\lambda)}
\left(\cos\pi(\lambda-u)A^+\left(u,v,\lambda\right)+
\cos\pi v A^-\left(u,v,\lambda\right)
\right),\label{S.alpha.compos}
\end{gather}
где
\begin{gather}\notag
C^{\pm}\left(u,v,\lambda\right)=
\sum_{d=1}^{\infty}
\frac{\tau_u(d)\tau_v(dN\pm1)\alpha(dN\pm1)}{d^{u}(dN\pm1)^{\lambda-u}}\Biggl(H_{\lambda}\left(u,v;\frac{\pm1}{dN\pm1}\right)-\\
-\frac{\Gamma(\lambda-u+v)\Gamma(\lambda-u-v)}{\Gamma(2\lambda)}\Biggr),\label{C.def}
\end{gather}
\begin{gather*}\label{A.def}
A^{\pm}\left(u,v,\lambda\right)=
\sum_{d=1}^{\infty}
\frac{\tau_u(d)\tau_v(dN\pm1)\alpha(dN\pm1)}{d^{u}(dN\pm1)^{\lambda-u}}.
\end{gather*}
Ввиду \eqref{H1.asym.for} ряды $C^{\pm}\left(u,v,\lambda\right)$ абсолютно сходятся при $\Re\lambda>0$ и, следовательно, аналитически продолжаются в точку $\lambda=1.$ Поэтому нам осталось продолжить в точку $\lambda=1$ только слагаемое
{\sloppy

}
\begin{gather}\label{A++A-}
\left(\cos\pi(\lambda-u)A^+\left(u,v,\lambda\right)+
\cos\pi v A^-\left(u,v,\lambda\right)
\right).
\end{gather}
Расписывая функцию $\tau_u(d)$ по определению \eqref{tau} и используя функцию $\eta(\cdot),$ определенную в \eqref{eta.def1}, получим, что
{\sloppy

}
\begin{gather}\notag
A^{\pm}\left(u,v,\lambda\right)=
\sum_{d_1,d_2=1}^{\infty}
\frac{\tau_v(d_1d_2N\pm1)\alpha(d_1d_2N\pm1)}{d_1^{2u}(d_1d_2N\pm1)^{\lambda-u}}
\left(\eta\left(\frac{d_1}{d_2}\right)+\eta\left(\frac{d_2}{d_1}\right)\right)=\\=
A_1^{\pm}\left(u,v,\lambda\right)+A_2^{\pm}\left(u,v,\lambda\right),\label{A=A1+A2}
\end{gather}
где первое слагаемое соответствует $\eta\left(d_1/d_2\right),$ а второе $\eta\left(d_2/d_1\right).$
В сумме $A_1^{+}\left(u,v,\lambda\right)$ сделаем внешним суммирование по $d_1,$ а в сумме по $d_2$ сделаем замену $m=d_1d_2N+1.$ В сумме $A_2^{+}\left(u,v,\lambda\right)$ сделаем внешним суммирование по $d_2,$ а в сумме по $d_1$ сделаем замену $m=d_1d_2N+1.$ В итоге получаем
\begin{gather*}\label{A1.transform1}
A_1^{+}\left(u,v,\lambda\right)=
\sum_{d_1=1}^{\infty}\frac{1}{d_1^{2u}}\sum_{m\equiv1(mod\, d_1N)}
\frac{\tau_v(m)\alpha(m)}{m^{\lambda-u}}
\eta\left(\frac{d_1^2N}{m-1}\right),
\end{gather*}
\begin{gather*}\label{A2.transform1}
A_2^{+}\left(u,v,\lambda\right)=N^{2u}
\sum_{d_2=1}^{\infty}d_2^{2u}\sum_{m\equiv1(mod\, d_2N)}
\frac{\tau_v(m)\alpha(m)}{m^{\lambda-u}(m-1)^{2u}}
\eta\left(\frac{d_2^2N}{m-1}\right).
\end{gather*}
Без ограничения общности можно считать, что $|u|\ll1,$ с абсолютной константой в знаке Виноградова. Так как
\begin{gather*}
(m-1)^{-2u}=m^{-2u}+O(1/m),
\end{gather*}
то
\begin{gather*}\label{A2.transform2}
A_2^{+}\left(u,v,\lambda\right)=N^{2u}
\sum_{d=1}^{\infty}d^{2u}\sum_{m\equiv1(mod\, dN)}
\frac{\tau_v(m)\alpha(m)}{m^{\lambda+u}}
\eta\left(\frac{d^2N}{m-1}\right)+O_{\varepsilon}\left(\frac{(PN)^{\varepsilon}}{PN}\right).
\end{gather*}
Заметим, что
\begin{gather*}
\eta(x)\neq\eta\left(x+\frac{x}{m-1}\right)\quad \hbox{только при}\quad
\frac{1}{2}-\frac{1}{2m}\le x\le 2
\end{gather*}
и
\begin{gather*}
\left|\eta(x)-\eta\left(x+\frac{x}{m-1}\right)\right|\ll\frac{x}{m-1}.
\end{gather*}
Так как
\begin{gather*}
\frac{d^2N}{m-1}=\frac{d^2N}{m}\left(1+\frac{1}{m-1}\right)
\end{gather*}
то
\begin{gather*}\label{A1.transform3}
A_1^{+}\left(u,v,\lambda\right)=
\sum_{d=1}^{\infty}\frac{1}{d^{2u}}\sum_{m\equiv1(mod\, dN)}
\frac{\tau_v(m)\alpha(m)}{m^{\lambda-u}}
\eta\left(\frac{d^2N}{m}\right)+O_{\varepsilon}\left(\frac{(PN)^{\varepsilon}}{PN}\right),
\end{gather*}
\begin{gather*}\label{A2.transform3}
A_2^{+}\left(u,v,\lambda\right)=N^{2u}
\sum_{d=1}^{\infty}d^{2u}\sum_{m\equiv1(mod\, dN)}
\frac{\tau_v(m)\alpha(m)}{m^{\lambda+u}}
\eta\left(\frac{d^2N}{m}\right)+O_{\varepsilon}\left(\frac{(PN)^{\varepsilon}}{PN}\right).
\end{gather*}
Аналогичные формулы с заменой $m\equiv1(mod\, dN)$ на $m\equiv-1(mod\, dN)$ справедливы и для $A_{1,2}^{-}\left(u,v,\lambda\right).$ Следовательно,
\begin{gather}\label{A2=A1}
A_2^{\pm}\left(u,v,\lambda\right)=N^{2u}A_1^{\pm}\left(-u,v,\lambda\right).
\end{gather}
Теперь к главным членам из $A_{1}^{\pm}\left(u,v,\lambda\right)$ применим формулу суммирования Вороного \eqref{Voronoi} с $M=1$ и
\begin{gather*}\label{f+def}
f^{\pm}(x)=\frac{\alpha(\pm x)}{({\pm} x)^{\lambda-u}}\eta\left(\frac{d^2N}{\pm x}\right).
\end{gather*}
В условиях формулы суммирования Вороного требовалось, чтобы  функция имела компактный носитель, отделенный от нуля. В нашем случае этого можно добиться, если представить $\alpha(x)$ в виде
{\sloppy

}
\begin{gather*}
\alpha(x)=\sum_{j>J}a_j(x), \quad\hbox{где}\quad \textit{supp}\, a_j(x)=[2^{j-1},2^{j+1}]
\end{gather*}
и сначала применить формулу суммирования Вороного с функциями $a_j(x),$ а потом уже просуммировать по $j.$ Для простоты изложения мы опустим этот технический момент и применим формулу суммирования Вороного с функцией $\alpha(x).$ В результате получим
\begin{gather}\label{A1Voronoi}
A_1^{\pm}\left(u,v,\lambda\right)=L\left(u,v,\lambda\right)+
L\left(u,-v,\lambda\right)+R^{\pm}\left(u,v,\lambda\right)+
O_{\varepsilon}\left(\frac{(PN)^{\varepsilon}}{PN}\right),
\end{gather}
где
\begin{gather}\label{Ldef}
L\left(u,v,\lambda\right)=
\frac{\zeta(1+2v)}{N}
\sum_{d=1}^{\infty}\frac{\gamma(2v;1;dN)}{d^{1+2u}}
\int_{0}^{\infty}\frac{\alpha(x)}{x^{\lambda-u-v}}
\eta\left(\frac{d^2N}{x}\right)dx,
\end{gather}
\begin{gather}\notag
R^{\pm}\left(u,v,\lambda\right)=
\frac{2\pi}{N^2}\sum_{q=1}^{\infty}\frac{1}{q^{2+2u}}
\sum_{m,n=1}^{\infty}\left(\frac{m}{n}\right)^{v}\\
\times\left(
g^{+}\left(16\pi^2\frac{mn}{(qN)^2} \right)S_{qN}(m,n;\pm1)+
g^{-}\left(16\pi^2\frac{mn}{(qN)^2} \right)S_{qN}(m,n;\mp1)\right)\label{R+-def}
\end{gather}
и
\begin{gather}\label{g+-Voronoi}
g^{\pm}\left(y\right)=\int_{0}^{\infty}
\kk^{\pm}(v;\sqrt{xy})f^{+}(x)dx.
\end{gather}
В формуле \eqref{Ldef} представим функцию $\eta\left(\frac{d^2N}{x}\right)$ с помощью обратного преобразования Меллина \eqref{Mellin.invers.def}
\begin{gather*}
\eta\left(\frac{d^2N}{x}\right)=\frac{1}{2\pi i}\int\limits_{\Re s=c}
\hat{\eta}(s)\frac{x^s}{d^{2s}N^s}ds,\quad \hbox{где}\quad c>0.
\end{gather*}
Ввиду абсолютной сходимости  интегралов и рядов в \eqref{Ldef} получаем
\begin{gather*}\label{Ltransform1}
L\left(u,v,\lambda\right)=
\frac{\zeta(1+2v)}{N}
\int_{0}^{\infty}\frac{\alpha(x)}{x^{\lambda-u-v}}
\frac{1}{2\pi i}\int\limits_{\Re s=c}
\hat{\eta}(s)\frac{x^s}{N^s}
\sum_{d=1}^{\infty}\frac{\gamma(2v;1;dN)}{d^{1+2u+2s}}dsdx.
\end{gather*}
Применяя к сумме по $d$ формулу \eqref{sumofgamma}, находим
\begin{gather*}\label{Ltransform2}
L\left(u,v,\lambda\right)=
\frac{\zeta(1+2v)}{N}
\int_{0}^{\infty}\frac{\alpha(x)}{x^{\lambda-u-v}}
\frac{1}{2\pi i}\int\limits_{\Re s=c}
\hat{\eta}(s)\frac{x^s}{N^s}
\frac{\zeta(1+2u+2s)}{\zeta(2+2u+2v+2s)}
\frac{\phi_{2v}(N)}{\phi_{1+2u+2v+2s}(N)}dsdx.
\end{gather*}
В интеграле по $s$ сдвинем контур на прямую $\sigma=-1/2.$ При этом мы пройдем полюса в точках $s_1=0$ и $s_2=u.$ Следовательно,
\begin{gather}\notag
L\left(u,v,\lambda\right)=
\frac{\zeta(1+2u)\zeta(1+2v)}{N}Z_1(u,v,\lambda)+\\+
\frac{\hat{\eta}(-u)\zeta(1+2v)}{2N^{1-u}}Z_1(0,v,\lambda)+
\frac{\zeta(1+2v)}{N}Z_2(u,v,\lambda),\label{Ltransform3}
\end{gather}
где
\begin{gather}\label{Z1def}
Z_1(u,v,\lambda)=\frac{\phi_{2v}(N)}{\zeta(2+2u+2v)\phi_{1+2u+2v}(N)}
\int_{0}^{\infty}\frac{\alpha(x)}{x^{\lambda-u-v}}dx,
\end{gather}
\begin{gather}\label{Z2def}
Z_2(u,v,\lambda)=
\int_{0}^{\infty}\frac{\alpha(x)}{x^{\lambda-u-v}}
\frac{1}{2\pi i}\int\limits_{\Re s=\sigma}
\hat{\eta}(s)\frac{x^s}{N^s}
\frac{\zeta(1+2u+2s)}{\zeta(2+2u+2v+2s)}
\frac{\phi_{2v}(N)}{\phi_{1+2u+2v+2s}(N)}dsdx.
\end{gather}
После интегрирования \eqref{Z1def} по частям получаем, ЧТО
\begin{gather}\label{Z1def2}
Z_1(u,v,\lambda)=\frac{\phi_{2v}(N)(\lambda-1-u-v)^{-1}}{\zeta(2+2u+2v)\phi_{1+2u+2v}(N)}
\int_{0}^{\infty}\frac{\alpha'(x)}{x^{\lambda-1-u-v}}dx.
\end{gather}
Напомним, что наша цель -- продолжить в точку $\lambda=1$ аналитическую функцию \eqref{A++A-}. Из \eqref{A=A1+A2}, \eqref{A2=A1} и \eqref{A1Voronoi} следует, что
\begin{gather*}\notag\label{A=L+L+L+L+R+R+O}
A^{\pm}\left(u,v,\lambda\right)=
L\left(u,v,\lambda\right)+L\left(u,-v,\lambda\right)+
N^{2u}L\left(-u,v,\lambda\right)+N^{2u}L\left(-u,-v,\lambda\right)+\\+
R^{\pm}\left(u,v,\lambda\right)+N^{2u}R^{\pm}\left(-u,v,\lambda\right)+
O_{\varepsilon}\left(\frac{(PN)^{\varepsilon}}{PN}\right).
\end{gather*}
Применяя \eqref{Ltransform3} и учитывая \eqref{etaMell+etaMell=0}, получаем
\begin{gather}\label{A=Abeaut+Rbeaut}
A^{\pm}\left(u,v,\lambda\right)=\A\left(u,v,\lambda\right)+\R^{\pm}\left(u,v,\lambda\right)
\end{gather}
\begin{gather}\label{Rbeaut=R+R+O}
\R^{\pm}\left(u,v,\lambda\right)=
R^{\pm}\left(u,v,\lambda\right)+N^{2u}R^{\pm}\left(-u,v,\lambda\right)+
O_{\varepsilon}\left(\frac{(PN)^{\varepsilon}}{PN}\right),
\end{gather}
\begin{gather}\notag
\A\left(u,v,\lambda\right)=
\frac{\zeta(1+2u)\zeta(1+2v)}{N}Z_1(u,v,\lambda)+
N^{2u}\frac{\zeta(1-2u)\zeta(1+2v)}{N}Z_1(-u,v,\lambda)+\\\notag+
\frac{\zeta(1+2u)\zeta(1-2v)}{N}Z_1(u,-v,\lambda)+
N^{2u}\frac{\zeta(1-2u)\zeta(1-2v)}{N}Z_1(-u,-v,\lambda)+\\\notag+
\frac{\zeta(1+2v)}{N}Z_2(u,v,\lambda)+
\frac{\zeta(1-2v)}{N}Z_2(u,-v,\lambda)+\\+
N^{2u}\frac{\zeta(1+2v)}{N}Z_2(-u,v,\lambda)+
N^{2u}\frac{\zeta(1-2v)}{N}Z_2(-u,-v,\lambda).\label{Abeaut=4Z1+4Z2}
\end{gather}
Из \eqref{R+-def} , \eqref{Z1def2} и \eqref{Z2def} следует, что при $u\neq0,\,v\neq0$ и $u+v\neq0$ мы можем продолжить $A^{\pm}\left(u,v,\lambda\right)$ в точку $\lambda=1.$

\begin{Le}\label{continuation lambda=1}
Для $k=1,\, \Re u=0,\,\Re v=0$ и $u\neq0,\,v\neq0$ и $u+v\neq0$ справедливо следующее тождество
\begin{gather}\label{continlambda=1}
\frac{\Gamma(2k-1)}{(4\pi)^{2k-1}}\sum_{f\in O_{2k}(N)}
L_f(\frac{1}{2}+u+v)L^{*}_f(\frac{1}{2}+u-v)=\zeta(1+2u)+\\\notag+
\frac{2(2\pi)^{2u}}{N^{2u}}\zeta(2u)\Gamma(2u)\cos\pi u
\frac{\Gamma(k-u+v)\Gamma(k-u-v)}{\Gamma(k+u+v)\Gamma(k+u-v)}+2\pi(-1)^k\delta_{1,N}d(u,v,k)
+\\\notag+
2\pi(-1)^kS_{\beta}(u,v;k)+
(-1)^k\frac{2(2\pi)^{2u}}{N^{2u}}
\left(
\cos\pi(k-u)C^+\left(u,v,k\right)+
\cos\pi v C^-\left(u,v,k;\right)\right)+\\\notag+
(-1)^k\frac{4(2\pi)^{2u}}{N^{2u}}
\frac{\Gamma(k-u+v)\Gamma(k-u-v)}{\Gamma(2k)}
\sin\frac{\pi(u+v)}{2}\sin\frac{\pi(u-v)}{2}\A\left(u,v,k\right)+\\\notag+
(-1)^k\frac{2(2\pi)^{2u}}{N^{2u}}
\frac{\Gamma(k-u+v)\Gamma(k-u-v)}{\Gamma(2k)}
\left(\cos\pi(k-u)\R^+\left(u,v,k\right)+
\cos\pi v \R^-\left(u,v,k\right)\right),
\end{gather}
где $d(u,v,k),\,S_{\beta}(u,v;k),\,C^{\pm}\left(u,v,k\right),\,\A\left(u,v,k\right)$ и $\R^{\pm}\left(u,v,k\right)$ определены соответственно в \eqref{d.u.v.lam}, \eqref{S.alpha.beta.def}, \eqref{C.def},\eqref{Abeaut=4Z1+4Z2} и \eqref{Rbeaut=R+R+O}.
\begin{proof}
Утверждение непосредственно следует из \eqref{basicformula1}, \eqref{D0func.int.repr.4}, \eqref{S.alpha.compos} и \eqref{A=Abeaut+Rbeaut}.
\end{proof}
\end{Le}
\section{Аналитическое продолжение в точку $u=0$}
Наша цель, получить аналитическое продолжение функции из правой части \eqref{continlambda=1} в точку $u=0.$ Отметим, что все
слагаемые в формуле \eqref{continlambda=1}, за исключением только первых двух и, возможно, функции $\A\left(u,v,k\right),$
являются аналитическими функциями. Сначала покажем, что сумма двух первых слагаемых  задает аналитическую функцию. Из функционального уравнение для дзета-функции Римана \eqref{zeta.func.equat} следует, что
{\sloppy

}
\begin{gather*}
\zeta(1+2u)+
\frac{2(2\pi)^{2u}}{N^{2u}}\zeta(2u)\Gamma(2u)\cos\pi u
\frac{\Gamma(k-u+v)\Gamma(k-u-v)}{\Gamma(k+u+v)\Gamma(k+u-v)}=\\=
\zeta(1+2u)+
\frac{(4\pi^2)^{2u}}{N^{2u}}\frac{\Gamma(k-u+v)\Gamma(k-u-v)}{\Gamma(k+u+v)\Gamma(k+u-v)}\zeta(1-2u).
\end{gather*}
Применяя \eqref{zeta.Loran} и переходя к пределу по $u\rightarrow0,$ получаем
\begin{gather*}\notag
\lim\limits_{u\rightarrow0}\left(
\zeta(1+2u)+
\frac{(4\pi^2)^{2u}}{N^{2u}}\frac{\Gamma(k-u+v)\Gamma(k-u-v)}{\Gamma(k+u+v)\Gamma(k+u-v)}\zeta(1-2u)\right)=\\=
\log N+2\gamma-2\log(2\pi)+\frac{\Gamma'(k+v)}{\Gamma(k+v)}+
\frac{\Gamma'(k-v)}{\Gamma(k-v)}=W_{2k}(N;iv),\label{u=0contin1}
\end{gather*}
где $W_{2k}(N;t)$ определено в \eqref{W}. Запишем $\A\left(u,v,k\right)$ в следующем виде
\begin{gather*}\label{Abeaut=Abeaut1+Abeaut2}
\A\left(u,v,k\right)=\A_1\left(u,v,k\right)+\A_1\left(u,-v,k\right)+\A_2\left(u,v,k\right)+\A_2\left(u,-v,k\right),
\end{gather*}
где
\begin{gather}\label{Abeaut1def}
\A_1\left(u,v,k\right)=
\frac{\zeta(1+2u)\zeta(1+2v)}{N}Z_1(u,v,\lambda)+
N^{2u}\frac{\zeta(1-2u)\zeta(1+2v)}{N}Z_1(-u,v,\lambda),
\end{gather}
\begin{gather}\label{Abeaut2}
\A_2\left(u,v,k\right)=
\frac{\zeta(1+2v)}{N}Z_2(u,v,\lambda)+
N^{2u}\frac{\zeta(1+2v)}{N}Z_2(-u,v,\lambda).
\end{gather}
Из \eqref{Z2def} следует, что $\A_2\left(u,v,k\right)$ аналитична в точке $u=0.$ Таким образом, нам осталось продолжить в точку $u=0$ только следующее выражение
\begin{gather}\label{h(A+A)}
\frac{4(-1)^k}{\Gamma(2k)}h(u,v)
\left(\A_1\left(u,v,k\right)+\A_1\left(u,-v,k\right)\right),
\end{gather}
где
\begin{gather*}\label{h.def}
h(u,v)=\left(\frac{2\pi}{N}\right)^{2u}\Gamma(k-u+v)\Gamma(k-u-v)\sin\frac{\pi(u+v)}{2}\sin\frac{\pi(u-v)}{2}.
\end{gather*}
Так как $h(u,v)=h(u,-v),$ то выражение \eqref{h(A+A)} симметрично по $v,$ и поэтому достаточно  продолжить в точку $u=0$ только слагаемое с $\A_1\left(u,v,k\right).$ Чтобы продолжить второе слагаемое с $\A_1\left(u,-v,k\right),$  достаточно заменить всюду $v$ на $-v.$ Из \eqref{Z1def2} \eqref{Abeaut1def} и  \eqref{h(A+A)} следует, что необходимо продолжить в точку $u=0$ выражение
\begin{gather*}
\frac{4(-1)^k}{\Gamma(2k)}h(u,v)
\A_1\left(u,v,k\right)=-\frac{4(-1)^k}{\Gamma(2k)}\frac{\zeta(1+2v)\phi_{2v}(N)}{N}h(u,v)\\
\int_{0}^{\infty}\alpha'(x)x^{v}
\left(\zeta(1+2u)\omega_1(u,v;x)+\zeta(1-2u)\omega_1(-u,v;x)N^{2u} \right)dx,
\end{gather*}
где
\begin{gather*}\label{omega1.def}
\omega_1(u,v;x)=\frac{x^u}{(u+v)\zeta(2+2u+2v)\phi_{1+2u+2v}(N)}.
\end{gather*}
Применяя \eqref{zeta.Loran} и переходя к пределу по $u\rightarrow0,$ получаем
\begin{gather*}
\lim\limits_{u\rightarrow0}h(u,v)
\left(\zeta(1+2u)\omega_1(u,v;x)+\zeta(1-2u)\omega_1(-u,v;x)N^{2u} \right)=\\=
(2\gamma-\log N)h(0,v)\omega_1(0,v;x)+h(0,v)\frac{\partial}{\partial u}\omega_1(u,v;x)\Bigl|_{u=0}
\end{gather*}
В итоге
\begin{gather}\label{Tbeaut.def}
\lim\limits_{u\rightarrow0}\frac{4(-1)^k}{\Gamma(2k)}h(u,v)\A_1\left(u,v,k\right)=
\frac{4(-1)^k}{\Gamma(2k)}h(0,v)\frac{\zeta(1+2v)\phi_{2v}(N)}{N}\\\notag
\times\int_{0}^{\infty}\alpha'(x)x^{v}
\left(
(\log N-2\gamma)\omega_1(0,v;x)-\frac{\partial}{\partial u}\omega_1(u,v;x)\Bigl|_{u=0}
\right)dx:=\frac{4(-1)^k}{\Gamma(2k)}h(0,v)\T(v,k).
\end{gather}
Итак, мы доказали следующий результат
\begin{Th}\label{continuation u=0}
Для $k=1,\,\Re v=0$ и $v\neq0$ справедлива формула
\begin{gather}\label{contin.u=0}
\frac{\Gamma(2k-1)}{(4\pi)^{2k-1}}\sum_{f\in O_{2k}(N)}
L_f(\frac{1}{2}+v)L^{*}_f(\frac{1}{2}-v)=W_{2k}(N;iv)+
2\pi(-1)^k\delta_{1,N}d(0,v,k)+\\\notag
+
2\pi(-1)^kS_{\beta}(0,v;k)+
2(-1)^k
\left(
(-1)^kC^+\left(0,v,k\right)+
\cos\pi v C^-\left(0,v,k;\right)\right)+\\\notag+
\frac{4(-1)^k}{\Gamma(2k)}h(0,v)
\left(\T(v,k)+\T(-v,k)+
\A_2\left(0,v,k\right)+\A_2\left(0,-v,k\right)\right)+\\\notag+
2(-1)^k
\frac{\Gamma(k+v)\Gamma(k-v)}{\Gamma(2k)}
\left((-1)^k\R^+\left(0,v,k\right)+
\cos\pi v \R^-\left(0,v,k\right)\right),
\end{gather}
где $W_{2k}(N;t),\,d(u,v,k),\,S_{\beta}(u,v;k),\,C^{\pm}\left(u,v,k\right),\,\R^{\pm}\left(u,v,k\right)$ и $\A_2\left(u,v,k\right)$ определены в \eqref{W}, \eqref{d.u.v.lam}, \eqref{S.alpha.beta.def}, \eqref{C.def},\eqref{Rbeaut=R+R+O} и \eqref{Abeaut2} соответственно.
\end{Th}
\section{Аналитическое продолжение в точку $v=0$}
Наша цель, получить аналитическое продолжение формулы \eqref{contin.u=0} в точку $v=0.$ Отметим, что только функции $d(u,v,k), \A_2\left(u,v,k\right)$ и $\T(v,k),$ возможно, не являются аналитическими функциями.  Сначала продолжим в точку $v=0$ функцию $d(u,v,k)$. Применяя \eqref{zeta.Loran} и переходя к пределу по $u\rightarrow0,$ получаем
{\sloppy

}
\begin{gather*}\label{d.u.v.continuat.v=0}
\lim\limits_{v\rightarrow0}
d(0,v,k)=\frac{1}{2\pi}W_{2k}(1;0).
\end{gather*}
Теперь продолжим в точку $v=0$ слагаемые с функцией $\A_2\left(0,v,k\right).$  Отметим,что так как $\Re v=0,$ то
\begin{gather}\label{h(0,v)value}
h(0,v)=-\frac{\pi}{2}v\tan\frac{\pi v}{2}.
\end{gather}
Таким образом, $h(0,v)$ имеет ноль второго порядка в точке $v=0$. Однако, из \eqref{Abeaut2} и \eqref{Z2def} следует, что
$\A_2\left(0,v,k\right)$ имеет полюс первого порядка  в точке $v=0.$ Поэтому
\begin{gather*}\label{Abeaut2.contin.v=0}
\lim\limits_{v\rightarrow0}h(0,v)\A_2\left(0,v,k\right)=0.
\end{gather*}
Нам осталось продолжить в точку $v=0$ слагаемые с функцией $\T(v,k),$ определенной в \eqref{Tbeaut.def}. Вычисляя производную функции $\omega_1(u,v;x)$ по $u$ в точке $u=0$ и учитывая \eqref{h(0,v)value}, получаем
\begin{gather}\label{Tbeaut.transform1}
\frac{4(-1)^k}{\Gamma(2k)}h(0,v)\T(v,k)=-\frac{2\pi(-1)^k}{N\Gamma(2k)}
\int_{0}^{\infty}\alpha'(x)\tan\frac{\pi v}{2}\zeta(1+2v)\left(\omega_2(v;x)+
\frac{1}{v}\omega_3(v;x)
\right)dx,
\end{gather}
где
\begin{gather}\label{omega2.def}
\omega_2(v;x)=\omega_3(v;x)
\left(
\log N-2\gamma-\log x+2\frac{\zeta'(2+2v)}{\zeta(2+2v)}+2\frac{\phi'_{1+2v}(N)}{\phi_{1+2v}(N)},
\right)
\end{gather}
\begin{gather}\label{omega3.def}
\omega_3(v;x)=
\frac{x^v\phi_{2v}(N)}{\zeta(2+2v)\phi_{1+2v}(N)}
\end{gather}
и $\phi'_{1+2v}(N)=\frac{\partial}{\partial s}\phi_s(N)|_{s=1+2v}.$  Вычисляя пределы по правилу Лопиталя, находим
\begin{gather*}\notag
\frac{4(-1)^k}{\Gamma(2k)}\lim\limits_{v\rightarrow0}
h(0,v)\left(\T(v,k)+\T(-v,k)\right)=\\
=-\frac{\pi^2(-1)^k}{N\Gamma(2k)}
\int_{0}^{\infty}\alpha'(x)\left(\omega_2(0;x)+2\gamma\omega_3(0;x)+
\frac{\partial}{\partial v}\omega_3(v;x)\Bigl|_{v=0}
\right)dx.\label{Tbeaut.limit v=0}
\end{gather*}
Используя \eqref{omega2.def} и \eqref{omega3.def}, получаем
\begin{gather*}\label{Tbeaut.limit2 v=0}
\frac{4(-1)^k}{\Gamma(2k)}\lim\limits_{v\rightarrow0}
h(0,v)\left(\T(v,k)+\T(-v,k)\right)=
-\frac{\pi^2(-1)^k}{N\Gamma(2k)}
\frac{\phi_{0}(N)}{\zeta(2)\phi_{1}(N)}
\left(
\log N+2\frac{\phi'_{0}(N)}{\phi_{0}(N)}
\right).
\end{gather*}
Таким образом, нами доказан следующий результат
\begin{Th}\label{continuation v=0}
Для $k=1$  справедлива формула
\begin{gather*}\notag
\frac{\Gamma(2k-1)}{(4\pi)^{2k-1}}\sum_{f\in O_{2k}(N)}
L_f(\frac{1}{2})L^{*}_f(\frac{1}{2})=W_{2k}(N;0)+
(-1)^k\delta_{1,N}W_{2k}(N;0)+\\\notag
+
2\pi(-1)^kS_{\beta}(0,0;k)+
2(-1)^k
\left(
(-1)^kC^+\left(0,0,k\right)+
 C^-\left(0,0,k;\right)\right)-\\\notag-
\frac{\pi^2(-1)^k}{\zeta(2)\Gamma(2k)}
\frac{\phi_{0}(N)}{N\phi_{1}(N)}
\left(
\log N+2\frac{\phi'_{0}(N)}{\phi_{0}(N)}
\right)+\\+
2(-1)^k
\frac{\Gamma(k)\Gamma(k)}{\Gamma(2k)}
\left((-1)^k\R^+\left(0,0,k\right)+
\R^-\left(0,0,k\right)\right),\label{contin.v=0}
\end{gather*}
где $W_{2k}(N;t),\,S_{\beta}(u,v;k),\,C^{\pm}\left(u,v,k\right)$ и $\R^{\pm}\left(u,v,k\right)$ определены в \eqref{W}, \eqref{S.alpha.beta.def}, \eqref{C.def} и \eqref{Rbeaut=R+R+O} соответственно.
\end{Th}
\begin{Zam}
Отметим, что при $N=1$ и нечетном $k$
\begin{gather*}
W_{2k}(N;0)+(-1)^k\delta_{1,N}W_{2k}(N;0)=0.
\end{gather*}
Это объясняется тем, что при нечетном $k$ из функционального уравнения \cite[лемма 3.6]{Mot.book} следует, что $L_f(\frac{1}{2})=0.$
{\sloppy

}
\end{Zam}
\section{Доказательство теоремы \ref{BFtheorem1}}
Для завершения доказательства теоремы \ref{BFtheorem1} нам осталось оценить слагаемые (все, кроме первых двух) из формулы \eqref{contin.u=0}. Положим $v=it.$ Чтобы оценить слагаемые $S_{\beta}(0,v;k),\,C^{\pm}\left(0,v,k\right),$ достаточно воспользоваться оценками на гипергеометрическую функцию Гаусса, доказанными в \cite{BF}. Получаем, что
\begin{gather}\label{SbetaC+-bound}
\left|S_{\beta}(0,v;k)\right|\ll_{\varepsilon}\frac{(1+|t|)P^{\varepsilon}}{N},\qquad
\left|C^+\left(0,v,k\right)+\cos\pi v C^-\left(0,v,k;\right)\right|\ll_{\varepsilon}\frac{(1+|t|)^3P^{\varepsilon}}{NP}
\end{gather}
Чтобы оценить слагаемое $h(0,v)\T(v,k),$ необходимо воспользоваться представлением \eqref{Tbeaut.transform1}. Так как
\begin{gather*}
|\zeta(1+2v)|\ll\log t,\qquad |\phi_{2v}(N)|\ll_{\varepsilon}N^{\varepsilon},
\end{gather*}
\begin{gather*}
|\phi_{1+2v}(N)|=\left|\prod_{p|N}\left(1-\frac{1}{p^{2+2v}}\right)\right|\ge
\prod_{p|N}\left(1-\frac{1}{p^{2}}\right)\gg1,
\end{gather*}
то
\begin{gather}\label{Tbeaut.bound}
\left|h(0,v)\T(v,k)\right|\ll_{\varepsilon}\frac{P^{\varepsilon}}{N}.
\end{gather}
Чтобы оценить слагаемое $h(0,v)\A_2(0,v,k),$ необходимо воспользоваться представлением \eqref{Abeaut2} и формулой \eqref{Z2def}. Используя \eqref{etaMellbound}, имеем
\begin{gather*}
\left|Z_2(0,v,k)\right|\ll_{\varepsilon} N^{1/2+\varepsilon}
\int_{0}^{\infty}\frac{\alpha(x)}{x^{k+1/2}}
\ll_{\varepsilon} \frac{N^{1/2+\varepsilon}}{P^{1/2}}.
\end{gather*}
Следовательно, применяя \eqref{h(0,v)value}, получаем
\begin{gather}\label{Abeaut_2.bound}
\left|h(0,v)\A_2(0,v,k)\right|\ll_{\varepsilon}\frac{(1+|t|)P^{\varepsilon}}{\left(PN\right)^{1/2}}.
\end{gather}
Таким образом нам осталось оценить в формуле \eqref{contin.u=0} только слагаемые с $\R^{\pm}\left(0,v,k\right).$ Используя формулу Стирлинга и \eqref{Rbeaut=R+R+O}, находим, что необходимо оценить величины
\begin{gather*}
(1+|t|)e^{-\pi|t|}R^+\left(0,v,k\right), \qquad (1+|t|)R^-\left(0,v,k\right).
\end{gather*}
Нам достаточно рассмотреть второе выражение, так как ввиду
формулы \eqref{R+-def} и  неравенства \eqref{Kloosterman.bound} выражение $(1+|t|)R^+\left(0,v,k\right)$ оценивается аналогично. Применяя \eqref{Kloosterman.bound}, получаем
{\sloppy

}
\begin{gather}\notag
(1+|t|)\left|R^{-}\left(0,v,k\right)\right|\ll_{\varepsilon}
\sum_{q=1}^{\infty}\frac{(1+|t|)}{(qN)^{3/2}}
\sum_{d=1}^{\infty} \tau(d)(d,qN)^{1/2}\\
\times\left(
g^{+}\left(16\pi^2\frac{d}{(qN)^2} \right)+
g^{-}\left(16\pi^2\frac{d}{(qN)^2} \right)\right),\label{R+-bound}
\end{gather}
где $g^{\pm}\left(y\right)$ определены в \eqref{g+-Voronoi}. Нам удобнее положить $v=it/2, t>0$ и
\begin{gather*}\label{theta.def}
\theta(x)=\alpha(x)\eta\left(\frac{d^2N}{x}\right).
\end{gather*}
Тогда
\begin{gather}\label{g+-Voronoi2}
g^{\pm}\left(y\right)=\int_{0}^{\infty}
\kk^{\pm}(v;\sqrt{xy})\theta(x)\frac{dx}{x}.
\end{gather}
Заметим, что $\theta(x)\neq0$ только при $x\ge\max\{P,q^2N/2\}=M.$ Оценим $g^{-}(y)$
\begin{Le}\label{lemma.g-.bound}
Положим $a=\frac{\sqrt{My}}{t}.$ Тогда
\begin{gather*}\label{g-bound.a>2}
|g^{-}(y)|\ll\frac{e^{-t\rho(a)}}{(at)^{3/2}},\qquad \hbox{если} \qquad a>2,
\end{gather*}
где $\rho(a)$ определена после \eqref{K.bound>1} и
\begin{gather*}\label{g-bound.a<2}
|g^{-}(y)|\ll
\frac{1+|\log a|}{t^{1/2}}\qquad \hbox{если} \qquad a\le2.
\end{gather*}
\begin{proof}
В интеграле \eqref{g+-Voronoi2} сделаем замену $xy=(tz)^2$ и воспользуемся \eqref{K+K-def}. В результате  получим
\begin{gather*}
|g^{-}(y)|\ll e^{\pi|t|/2}\left|\int_{0}^{\infty}
K_{it}(tz)\theta\left(\frac{(tz)^2}{y}\right)\frac{dz}{z}\right|.
\end{gather*}
Заметим, что $\theta\left(\frac{(tz)^2}{y}\right)\neq0$ только при $z\ge\frac{\sqrt{My}}{t}=a.$ Возможны три случая.
\begin{enumerate}
  \item Если $a>2,$ тогда, используя \eqref{K.bound>1} получаем
\begin{gather*}
|g^{-}(y)|\ll\int_{a}^{\infty}
e^{-t\rho(z)}\frac{1}{t^{1/2}z(z^2-1)^{1/4}}dz=
\frac{1}{t^{1/2}}\int_{a}^{\infty}
e^{-t\rho(z)}\rho'(z)\frac{dz}{(z^2-1)^{3/4}}.
\end{gather*}
Применяя вторую теорему о среднем, находим
\begin{gather*}
|g^{-}(y)|\ll\frac{e^{-t\rho(a)}}{(at)^{3/2}}.
\end{gather*}
  \item Если $1<a\le2,$ тогда, используя \eqref{K.bound>1} получаем
\begin{gather}\label{g-bound.1<a<2}
|g^{-}(y)|\ll\int_{1}^{1+t^{-2/3}}
e^{-t\rho(z)}\frac{dz}{zt^{1/3}}+
\int_{1+t^{-2/3}}^{\infty}
e^{-t\rho(z)}\frac{dz}{t^{1/2}z(z^2-1)^{1/4}}\ll\frac{1}{t}.
\end{gather}
  \item Если $a\le1,$ тогда, используя \eqref{K.bound<1} и \eqref{g-bound.1<a<2} находим
\begin{gather*}
|g^{-}(y)|\ll\int_{a}^{1}
\min\left\{\frac{1}{t^{1/3}},\frac{1}{t^{1/2}(1-z^2)^{1/4}} \right\}\frac{dz}{z}+\frac{1}{t}\ll
\frac{1+|\log a|}{t^{1/2}}.
\end{gather*}
\end{enumerate}
\end{proof}
\end{Le}
\begin{Le}\label{lemma.g+.bound}
Справедливы оценки
\begin{gather*}\label{g+bound}
|g^{+}(y)|\ll\frac{1}{t^{5/2}}
\left\{
              \begin{array}{ll}
               1+|\log a|, & \hbox{если $a<1$} \\
                (1+a^2)^{-5/4}, & \hbox{если $a>1$.}
              \end{array}
\right.
\end{gather*}
\begin{proof}
Из \cite[8.491, 4)]{GradRyz} следует, что
\begin{gather}\label{k+diff.equation}
\frac{\partial^2}{\partial x^2}\left(\sqrt{x}\kk^+(v;\sqrt{xy})\right)=
-\frac{xy+1+t^2}{4x^2}\left(\sqrt{x}\kk^+(v;\sqrt{xy})\right)
\end{gather}
Интегрируя \eqref{g+-Voronoi2} два раза по частям с помощью \eqref{k+diff.equation}, получаем
\begin{gather*}
|g^{+}(y)|\ll\left|\int_{0}^{\infty}
\sqrt{x}\kk^{\pm}(v;\sqrt{xy})
\frac{\partial^2}{\partial x^2}\left(
\frac{x^{1/2}}{xy+1+t^2}\theta(x)
\right)dx\right|.
\end{gather*}
Вычисляя вторую производную, получаем
\begin{gather*}\label{by part1}
|g^{+}(y)|\ll\left|\int_{0}^{\infty}
\kk^{\pm}(v;\sqrt{xy})
\sum_{n=0}^{2}
\frac{\theta^{(n)}(x)x^{n-1}}{xy+1+t^2}P_{2-n}\left(\frac{xy}{xy+1+t^2}\right)dx\right|,
\end{gather*}
где $P_k(z)$--многочлен степени $k.$ Отметим, что $\left|P_{2-n}\left(\frac{xy}{xy+1+t^2}\right)\right|\ll1.$
Из определения функций $\alpha(\cdot),\,\eta(\cdot),\,\theta(\cdot)$  следует, что при $n>0$
\begin{gather*}
\theta^{(n)}(x)x^{n}\asymp
\left\{
              \begin{array}{ll}
               1, & \hbox{если $M\le x\le2M$} \\
                0, & \hbox{иначе.}
              \end{array}
\right.
\end{gather*}
Таким образом, мы получаем
\begin{gather*}
|g^{+}(y)|\ll\int_{0}^{\infty}
|\kk^{\pm}(v;\sqrt{xy})|
\frac{\theta(x)}{x(xy+1+t^2)}dx.
\end{gather*}
Сделав замену $xy=(tz)^2$ и применяя оценку \eqref{K+bound}, находим
\begin{gather*}
|g^{+}(y)|\ll\frac{1}{t^{1/2}}\int_{0}^{\infty}
\frac{\theta\left(\frac{(tz)^2}{y}\right)}{(tz)^2+1+t^2}\frac{dz}{z(1+z^2)^{1/4}}
\end{gather*}
Так как $\theta\left(\frac{(tz)^2}{y}\right)\neq0$ только при $z\ge a,$ то
\begin{gather*}
|g^{+}(y)|\ll\frac{1}{t^{5/2}}
\left\{
              \begin{array}{ll}
               1+|\log a|, & \hbox{если $a<1$} \\
                (1+a^2)^{-5/4}, & \hbox{если $a>1$.}
              \end{array}
\right.
\end{gather*}

\end{proof}
\end{Le}
Подставляя в \eqref{R+-bound} оценки из лемм \ref{lemma.g-.bound} и \ref{lemma.g+.bound},  получаем
\begin{gather}\label{R+-bound2}
(1+|t|)\left|R^{-}\left(0,v,k\right)\right|\ll_{\varepsilon}
\frac{(1+|t|)^{5/2}}{(PN)^{1/4}}P^{\varepsilon}.
\end{gather}
Положив $P=(1+|t|)^{6}N^{3}$, из формулы \eqref{contin.u=0} и оценок \eqref{SbetaC+-bound}, \eqref{Tbeaut.bound}, \eqref{Abeaut_2.bound} и \eqref{R+-bound2} получаем теорему \ref{BFTH1}.

\section{Доказательство теоремы \ref{BFtheorem2}}
Из оценок \eqref{SbetaC+-bound}, \eqref{Tbeaut.bound}, \eqref{Abeaut_2.bound} и \eqref{R+-bound2}, а также выбора $P=(1+|t|)^{a}N^{b}$ с достаточно большими $a,b$ вытекает, что для доказательства теоремы \ref{BFtheorem2} нам достаточно получить оценку
\begin{gather*}\label{integral of Sbeta-bound}
\int_{0}^{T}\left|S_{\beta}(0,it;k)\right|dt\ll_{\varepsilon}\frac{T}{N}(TN)^{\varepsilon}.
\end{gather*}
В статье \cite{BF} была доказана оценка
\begin{gather*}
\left|H_{k}\left(0,v;\frac{1}{dN+1}\right)\right|\ll 1.
\end{gather*}
Следовательно,
\begin{gather*}
\int_{0}^{T}
\left|
\sum_{d=1}^{\infty}
\frac{\tau_0(d)\tau_v(dN+1)}{dN+1}H_{1}\left(0,v;\frac{1}{dN+1}\right)\beta(dN+1)\right|dt
\ll_{\varepsilon}\frac{T}{N}(TN)^{\varepsilon}.
\end{gather*}
Таким образом нам осталось доказать, что
\begin{gather}\label{integral of H- bound}
\left|\int_{T}^{2T}
\cos\pi v
\sum_{d=1}^{\infty}
\frac{\tau_0(d)\tau_v(dN-1)}{dN-1}H_{1}\left(0,v;\frac{-1}{dN-1}\right)\beta(dN-1)dt\right|
\ll_{\varepsilon}\frac{T}{N}(TN)^{\varepsilon}.
\end{gather}
Разобьем сумму по $d$ на две $\Sigma_1, \Sigma_2$ с условиями $dN-1\le 100T^2$ и $dN-1>100T^2$.  Если $100T^2<N,$ то первая сумма пуста. При работе с первой суммой мы расписываем $\tau_v(dN-1)$ по определению \eqref{tau} и используем интегральное представление для
$H_{k}\left(0,v;-x\right)$ из работы \cite[теорема 2]{BF2}
\begin{gather*}
\cos\pi v H_{k}\left(0,v;-x\right)=
\frac{2}{\sqrt{x}}\coth(\pi t)\int_{0}^{1}
\frac{\sin\left( (2k-1)\arcsin z\right)}{\sqrt{1-z^2}\sqrt{1+xz^2}}
\sin\left(2t \kappa\left(z\sqrt{x}\right)\right)dz,
\end{gather*}
где $\kappa(z)=\log(z+\sqrt{1+z^2}).$ В результате получаем
\begin{gather}\label{Sigma1def}
\Sigma_1=2
\sum_{dN-1\le T^2}
\frac{\tau_0(d)}{\sqrt{dN-1}}\sum_{mn=dN-1}
\int_{0}^{1}
\frac{\sin\left( (2k-1)\arcsin z\right)}{\sqrt{1-z^2}\sqrt{1+z^2/(dN-1)}}
I_T\left(z,\frac{m}{n},dN-1\right)dz,
\end{gather}
\begin{gather}\label{IT.def}
I_T\left(z,\frac{m}{n},dN-1\right)=\int_{T}^{2T}\coth(\pi t)
\left(\frac{m}{n}\right)^{it}
\sin\left(2t \kappa\left(\frac{z}{\sqrt{dN-1}}\right)\right)
dt.
\end{gather}
Для упрощения записи мы опускаем константу $100$ в суммировании по $d,$ так как она (константа) не оказывает влияния на ход вычислений.
Возможны несколько случаев:
\begin{enumerate}
  \item Если $|m-n|\ll 1,$ то существует точка $z_0\in[0,1]$ такая, что
\begin{gather}\label{station.phase1}
\log\frac{m}{n}\pm2\kappa\left(\frac{z_0}{\sqrt{dN-1}}\right)=0.
\end{gather}
Положим $\delta=\sqrt{dN-1}/T^{1+\epsilon}$ и разобьем интеграл по $z$ на два $I_1$ и $I_2$ с условиями  $|z-z_0|\le \delta$ и $|z-z_0|>\delta,$ соответственно. При $|z-z_0|\le \delta$ интеграл по \eqref{IT.def} оценим тривиально
\begin{gather*}
\left|I_T\left(z,\frac{m}{n},dN-1\right)\right|\ll T.
\end{gather*}
Так как $\kappa'\left(\frac{z}{\sqrt{dN-1}}\right)\asymp\frac{1}{\sqrt{dN-1}},$ то при $|z-z_0|>\delta$
\begin{gather*}
\left|I_T\left(z,\frac{m}{n},dN-1\right)\right|\ll \frac{\sqrt{dN-1}}{|z-z_0|}.
\end{gather*}
За счет выбора параметра $\delta$ мы получаем
\begin{gather}\label{I1+I2.bound}
I_1+I_2\ll T\delta+\sqrt{dN-1}\log\delta\ll \sqrt{dN-1}T^{\epsilon}.
\end{gather}
Подставляя оценку \eqref{I1+I2.bound} в \eqref{Sigma1def}, находим
\begin{gather*}
\Sigma_1\ll_{\epsilon} T^{\epsilon}
\sum_{dN-1\le T^2}\sum_{|c|\ll1}
\sum_{n(n+c)=dN-1}1=T^{\epsilon}
\sum_{|c|\ll1}
\sum\limits_{n(n+c)\equiv-1(mod N)\atop n(n+c)\le T^2}1.
\end{gather*}
Разбивая интервал суммирования по $n$ на отрезки вида $[jN,(j+1)N-1]$ и учитывая то, что квадратичное сравнение имеет конечное число решений, получаем
$
\Sigma_1\ll_{\epsilon}T^{1+\epsilon}/N.
$
  \item Если $1\ll|m-n|\ll \epsilon\sqrt{dN-1},$ то левая часть \eqref{station.phase1} не обнуляется, но остается достаточно малой. Для $c=m-n$
\begin{gather*}
\left|\log\frac{m}{n}\pm2\kappa\left(\frac{z}{\sqrt{dN-1}}\right)\right|\gg\left|\log\frac{m}{n}\right|
\gg\frac{c}{n}
\end{gather*}
и, следовательно,
\begin{gather*}
\left|I_T\left(z,\frac{m}{n},dN-1\right)\right|\ll \frac{n}{c}.
\end{gather*}
В результате находим
\begin{gather*}
\Sigma_1\ll_{\epsilon}
\sum_{dN-1\le T^2}\frac{T^{\epsilon}}{\sqrt{dN-1}}\sum_{1\ll|c|\ll\epsilon\sqrt{dN-1}}
\sum_{n(n+c)=dN-1}\frac{n}{c}\ll_{\epsilon}
\\
\ll_{\epsilon}T^{\epsilon}
\sum_{1\ll|c|\ll\epsilon T}\frac{1}{c}
\sum\limits_{n(n+c)\equiv-1(mod N)\atop n(n+c)\le T^2}1.\ll_{\epsilon}
\frac{T^{1+\epsilon}}{N}.
\end{gather*}
  \item
  Если $|m-n|\gg\epsilon\sqrt{dN-1},$ то левая часть \eqref{station.phase1} становится достаточно большой, и поэтому
$
\left|I_T\left(\cdot\right)\right|\ll 1.
$
Следовательно,
\begin{gather*}
\Sigma_1\ll_{\epsilon}
\sum_{dN-1\le T^2}\frac{T^{\epsilon}}{\sqrt{dN-1}}\ll_{\epsilon}
\frac{T^{1+\epsilon}}{N}.
\end{gather*}
\end{enumerate}
 При работе со второй суммой $\Sigma_2$ мы расписываем $\tau_v(dN-1)$ по определению \eqref{tau} и используем представление функции
$H_{1}\left(0,v;-x\right)$ в виде ряда
\begin{gather*}
H_{1}\left(0,v;-x\right)=
\sum_{l=0}^{\infty}\frac{\Gamma(l+1+v)\Gamma(l+1-v)}{\Gamma(l+1)\Gamma(l+2)}(-x)^l.
\end{gather*}
Данный ряд абсолютно сходится (см. \cite{BF}) при $(1+t)\sqrt{x}<1.$ Именно для того, чтобы это условие выполнялось, мы ввели ограничение $dN-1>100T^2.$ Таким образом,
\begin{gather*}\label{Sigma2def}
\Sigma_2=
\sum_{l=0}^{\infty}\frac{(-1)^l}{\Gamma(l+1)\Gamma(l+2)}
\sum_{dN-1>100T^2}
\frac{\tau_0(d)\beta(dN-1)}{(dN-1)^{l+1}}\sum_{mn=dN-1}
J_T\left(\frac{m}{n};l\right),
\end{gather*}
\begin{gather*}\label{JT.def}
J_T\left(\frac{m}{n};l\right)=\int_{T}^{2T}
\left(\frac{m}{n}\right)^{it}g(t)dt, \qquad
g(t)=\cosh(\pi t)\Gamma(l+1+it)\Gamma(l+1-it).
\end{gather*}
Отметим, что (см.\cite[5.4.3, 5.5.1]{NIST})
\begin{gather*}
\Gamma(l+1+it)\Gamma(l+1-it)=\left|\Gamma(it)\right|^2\prod_{j=0}^{l}(j^2+t^2)=
\frac{\pi}{t\sinh\pi t}\prod_{j=0}^{l}(j^2+t^2).
\end{gather*}
Положив $c=m-n,$ по второй теореме о среднем получаем
\begin{gather*}\label{Gamma*Gamma}
\left|J_T\left(\frac{m}{n};l\right)\right|\ll g(2T)
\left\{
              \begin{array}{ll}
               T,   & \hbox{если $c=0$} \\
               n/c, & \hbox{если $|c|\le \epsilon\sqrt{dN-1}$}\\
               1,   & \hbox{если $|c|> \epsilon\sqrt{dN-1}$}.
              \end{array}
\right.
\end{gather*}
Рассмотрим отдельно каждый из трех случаев
\begin{enumerate}
  \item Если $c=0,$ то
\begin{gather*}
\Sigma_2\ll
\sum_{l=0}^{\infty}\frac{Tg(2T)}{\Gamma(l+1)\Gamma(l+2)}
\sum\limits_{dN-1>100T^2\atop dN-1=n^2}
\frac{\tau_0(d)\beta(dN-1)}{(dN-1)^{l+1}}\ll_{\epsilon}\\\ll_{\epsilon}
\sum\limits_{100T^2<dN-1<P\atop dN-1=n^2}\frac{1}{(dN-1)^{1-\epsilon}}
\sum_{l=0}^{\infty}\frac{Tg(2T)}{\Gamma(l+1)\Gamma(l+2)}\frac{1}{(dN-1)^{l}}.
\end{gather*}
Используя оценки на функцию $g(T)$ из \cite{BF}, получаем
\begin{gather*}
\Sigma_2\ll_{\epsilon}T^2
\sum\limits_{100T^2<dN-1<P\atop dN-1=n^2}\frac{1}{(dN-1)^{1-\epsilon}}\ll_{\epsilon}T^2
\sum_{T/N\le j\le P/N}\sum\limits_{Nj<n\le N(j+1)\atop n^2\equiv -1(mod\,N)}\frac{1}{n^{2-\epsilon}}\ll_{\epsilon}\\
\ll_{\epsilon}T^2\sum_{T/N\le j\le P/N}\frac{1}{(Nj)^{2-\epsilon}}\ll_{\epsilon}\frac{T^{1+\epsilon}}{N}.
\end{gather*}
  \item Если $|c|\le \epsilon\sqrt{dN-1},$ то
\begin{gather*}
\Sigma_2\ll
\sum_{l=0}^{\infty}\frac{g(2T)}{\Gamma(l+1)\Gamma(l+2)}
\sum_{dN-1>100T^2}
\frac{\tau_0(d)\beta(dN-1)}{(dN-1)^{l+1}}\sum\limits_{0<|c|\le\epsilon\sqrt{dN-1}\atop dN-1=n(n+c)}\frac{n}{c}
\ll_{\epsilon}\\\ll_{\epsilon} T
\sum_{100T^2<dN-1<P}\frac{1}{(dN-1)^{1-\epsilon}}
\sum\limits_{0<|c|\le\epsilon\sqrt{dN-1}\atop dN-1=n(n+c)}\frac{n}{c}\ll_{\epsilon}\\
\ll_{\epsilon} T \sum_{0<|c|\le\epsilon\sqrt{P}}\frac{1}{c}
\sum\limits_{T<n<P^{1/2}\atop n(n+c)\equiv -1 (mod\,N)}\frac{1}{n^{1-\epsilon}}\ll_{\epsilon}\frac{TP^{\epsilon}}{N}.
\end{gather*}
  \item Если $|c|>\epsilon\sqrt{dN-1},$ то
\begin{gather*}
\Sigma_2\ll
\sum_{l=0}^{\infty}\frac{g(2T)}{\Gamma(l+1)\Gamma(l+2)}
\sum\limits_{dN-1>100T^2}
\frac{\tau_0(d)\beta(dN-1)}{(dN-1)^{l+1}}
\ll_{\epsilon}\\\ll_{\epsilon} T
\sum\limits_{100T^2<dN-1<P}\frac{1}{(dN-1)^{1-\epsilon}}\ll_{\epsilon}\frac{TP^{\epsilon}}{N}.
\end{gather*}
\end{enumerate}
Таким образом оценка \eqref{integral of H- bound} доказана, а, следовательно, доказана и теорема \ref{BFtheorem2}.



\noindent
\textit{V.A. Bykovskii}\\
\textit{Institute for Applied Mathematics, Russian Academy of Sciences, Khabarovsk, Russia,}\\
\textit{e-mail: vab@iam.khv.ru}\\
\textit{D.A. Frolenkov}\\
\textit{Steklov Mathematical Institute of Russian Academy of Sciences, Moscow, Russia,}\\
\textit{Institute for Applied Mathematics, Russian Academy of Sciences, Khabarovsk, Russia,}\\
\textit{e-mail: frolenkov@mi.ras.ru}

\end{document}